\documentclass[12pt,a4paper,fleqn]{article}
\usepackage[latin1]{inputenc}
\usepackage[T1]{fontenc}
\usepackage{amsmath}
\usepackage{amssymb}
\usepackage{textcomp}
\usepackage[english,frenchb]{babel}
\usepackage[all,dvips]{xy}

\newtheorem{thm}{Théorème}[section]
\newtheorem{lem}[thm]{Lemme}
\newtheorem{prop}[thm]{Proposition}
\newtheorem{cor}[thm]{Corollaire}

\newcommand\dem{\textbf{Démonstration: }}
\newcommand\Def{\textbf{Définition: }}
\newcommand\Defs{\textbf{Définitions: }}
\newcommand\rem{\textbf{Remarque: }}
\newtheorem{soulem}[thm]{Sous-lemme}

\title{Presque-réductibilité des cocycles quasi-périodiques de classe Gevrey 2}
\author{Claire Chavaudret\\
\small{Institut de Mathématiques de Jussieu}}
\date{}

\begin{document}

\maketitle

\section{Introduction}

\subsection{Présentation du résultat}

Soit $n\geq 1$, 
$d\geq 1,\omega\in\mathbb{R}^d$ 
et $0<\kappa<1,\tau>\max(1,d-1)$. 
Supposons que $\omega$ est diophantien de constante $\kappa$ et d'exposant $\tau$, 
c'est-à-dire que 

\begin{equation}
\forall \ m\in\mathbb{Z}^d\setminus\{0\},\ |\langle m,\omega\rangle|\geq \frac{\kappa}{|m|^\tau}\end{equation}

\noindent On peut aussi supposer sans perte de généralité 
que $\sup |\omega_i|\leq 1$. 
On note $\mathbb{T}^d:=\mathbb{R}^d/\mathbb{Z}^d$ le tore et $2\mathbb{T}^d:=\mathbb{R}^d/(2\mathbb{Z}^d)$ 
le double tore.

\bigskip
\noindent \Def Soit $\mathcal{G}$ une algèbre de Lie, $G$ le groupe de Lie associé à $\mathcal{G}$ et $A:2\mathbb{T}^d\rightarrow \mathcal{G}$. Le cocycle quasi-périodique 
associé à $A$ est la fonction $X:2\mathbb{T}^d\times \mathbb{R}\rightarrow G$ telle que pour tous 
$(\theta,t)\in 2\mathbb{T}^d\times \mathbb{R}$,

\begin{equation} \frac{d}{dt}X^t(\theta)=A(\theta+t\omega)X^t(\theta); \ X^0(\theta)=Id
\end{equation}

\bigskip
\noindent On dit que c'est un cocycle constant si $A$ est constante. 

\bigskip
\rem Le terme quasi-périodique vient du fait que $A$ est la fonction enveloppe d'une fonction quasi-périodique.
En effet, pour tout $\theta\in 2\mathbb{T}^d$, $t\mapsto A(\theta+t\omega ) $ est une fonction quasi-périodique.

\bigskip
\noindent Un cocycle constant est toujours de la forme $t\mapsto e^{tA}$. 

\bigskip
\noindent  
\Def Soit $r>0$. Soit $E$ un sous-ensemble de $gl(n,\mathbb{C})$. 
Notons $C^{G,2}_r(2\mathbb{T}^d,E)$ 
(resp. $C^{G,2}_r(\mathbb{T}^d,E)$) 
les fonctions de classe Gevrey deux de paramètre $r$ sur $2\mathbb{T}^d$ (resp. $\mathbb{T}^d$) à valeurs dans 
$E$, 
c'est-à-dire les fonctions $f\in
C^\infty(2\mathbb{T}^d,E)$ 
(resp. $C^\infty(\mathbb{T}^d,E)$) 
telles qu'il existe $C>0$ tel que pour tout $\alpha=(\alpha_1,\dots,\alpha_d)\in \mathbb{N}^d$,

\begin{equation}\sup_{\theta\in 2\mathbb{T}^d}|\partial^\alpha f(\theta)|\leq Cr^{-2\mid \alpha\mid }(\alpha!)^2
\end{equation}

\noindent où on a noté $\alpha!=\alpha_1!\dots \alpha_d!$. Notons $||f||_r$ la norme Gevrey de $f$, c'est-à-dire

\begin{equation}||f||_r=\sup_{\alpha\in \mathbb{N}^d} r^{2\mid \alpha\mid} \frac{1}{(\alpha!)^2} \sup_{\theta\in 2\mathbb{T}^d}| |\partial^\alpha f(\theta)| |
\end{equation}

\bigskip
\noindent \textbf{Notation:} Pour toute fonction $f\in C^1(2\mathbb{T}^d,\mathbb{C})$, on notera pour tout $\theta\in 2\mathbb{T}^d$

\begin{equation}\partial_\omega f(\theta)=\frac{d}{dt} f(\theta+t\omega)_{\mid t=0}
\end{equation}

\noindent 
la dérivée de $f$ dans la direction $\omega$.

\bigskip 
\noindent 
\Def Soit $\mathcal{G}$ une algèbre de Lie et $G$ le groupe de Lie associé à $\mathcal{G}$. Soient $r,r'>0$ et $A,B\in C^{G,2}_r(2\mathbb{T}^d,\mathcal{G})$. 
On dit que $A$ et $B$ sont \textit{conjugués dans $C^{G,2}_{r'}(2\mathbb{T}^d,G)$} s'il existe $Z\in C^{G,2}_{r'}(2\mathbb{T}^d,G)$ tel que pour 
tout $\theta\in 
2\mathbb{T}^d$,

$$\partial_\omega Z(\theta)=A(\theta)Z(\theta)-Z(\theta)B(\theta)$$

\noindent Si $B$ est constante en $\theta$, on dit que $A$ est \textit{réductible dans $C^{G,2}_{r'}(2\mathbb{T}^d,G)$},
ou \textit{réductible par $Z$ à 
$B$}.

\noindent 
\rem Soit $X$ le cocycle quasi-périodique associé à $A$. La fonction $A$ est réductible par $\Phi$ à $A_0$ si et seulement si  

\begin{equation}\forall (t,\theta),\ X^t(\theta)=\Phi(\theta+t\omega)^{-1}e^{tA_0}\Phi(\theta)
\end{equation}

\noindent La réductibilité équivaut aussi au fait que l'application de $2\mathbb{T}^d\times \mathbb{R}^n$ dans lui-même:

\begin{equation}
\left(\begin{array}{c}
\theta\\
v\\
\end{array}\right)
\mapsto
\left(\begin{array}{c}
\theta+\omega\\
X^1(\theta)v\\
\end{array}\right)
\end{equation}

\noindent soit conjuguée à une application $\chi$ telle que 

\begin{equation}\frac{d\chi}{d\theta}  \left(\begin{array}{c}
\theta\\
v\\
\end{array}\right)
\equiv \left(\begin{array}{c}
\bar{1}\\
{0}\\
\end{array}\right)\end{equation}

\bigskip
\noindent Le but de cet article est de montrer que pour 

$$G=GL(n,\mathbb{C}),GL(n,\mathbb{R}),SL(n,\mathbb{R}),Sp(n,\mathbb{R}),O(n),U(n) $$

\noindent au voisinage d'un cocycle constant, tout cocycle de classe Gevrey 2 de paramètre $r$ à valeurs dans $G$
 est 
presque-réductible dans $\cup_{r'>0}C^{G,2}_{r'}(2\mathbb{T}^d,G)$. 

\bigskip
\noindent Nous allons prouver le théorème suivant, pour $G$ parmi les groupes mentionnés ci-dessus et $\mathcal{G}$ l'algèbre de Lie associée à $G$:


\begin{thm}\label{th1}Soit $0<r\leq 1$, $A\in \mathcal{G}$, $F\in C^{G,2}_r(\mathbb{T}^d, \mathcal{G})$. 
Il existe $\epsilon_0<1$ 
ne dépendant que de $n,d,\kappa,\tau,A,r$ tel que si 

$$| | F| |_r\leq \epsilon_0$$ 

\noindent 
alors 
pour tout $\epsilon>0$, il existe $r_\epsilon>0,\bar{A}_\epsilon,\bar{F}_\epsilon
\in C^{G,2}_{r_\epsilon}
(2\mathbb{T}^d,\mathcal{G})$, 
$Z_\epsilon\in C^{G,2}_{r_\epsilon}(2\mathbb{T}^d, G)$ 
 tels que  pour tout $\theta\in 2\mathbb{T}^d$,

$$\partial_\omega Z_\epsilon(\theta)=(A+F(\theta))Z_\epsilon(\theta)-Z_\epsilon(\theta)(\bar{A}_\epsilon(\theta)
+\bar{F}_\epsilon(\theta))
$$

\noindent où 
\begin{itemize}
\item $\bar{A}_\epsilon$ est réductible dans $C^{G,2}_{r_\epsilon}(2\mathbb{T}^d, G)$, 
\item $| |\bar{F}_\epsilon| |_{r_\epsilon}\leq \epsilon$, 

\item et $| |Z_\epsilon
-Id | |_{r_\epsilon}\leq 2\epsilon_0^\frac{1}{2}$. 
\end{itemize}

\bigskip
\noindent De plus, en dimension 2 ou si $G=GL(n,\mathbb{C})$ ou $U(n)$, $Z_\epsilon,\bar{A}_\epsilon,\bar{F}_\epsilon$ 
sont continus sur $\mathbb{T}^d$. 

\end{thm}


\subsection{Généralisations et conséquences}

\noindent Le théorème \ref{th1} dit qu'au voisinage d'un 
cocycle constant, tous les cocycles sont presque-réductibles au sens où ils sont arbitrairement proches d'un cocycle réductible, ce qui signifie que la réductibilité est 
un phénomène prédominant. Cependant, si la réductibilité implique la presque-réductibilité, l'inverse n'est pas vrai: il existe des cocycles non réductibles même proches d'un cocycle constant. 

\bigskip
\noindent L'intérêt de la notion de presque-réductibilité est qu'un cocycle presque réductible a une dynamique connue sur un temps très long.

\bigskip
\noindent Le théorème \ref{th1} est en fait vrai si l'on suppose $F$ dans une classe plus grande que $C^{G,2}_r(\mathbb{T}^d,
\mathcal{G})$, 
à savoir les fonctions de $C^{G,2}_r(2\mathbb{T}^d,\mathcal{G})$ 
vérifiant certaines 
"bonnes propriétés de périodicité" par rapport à la matrice $A$. 

\bigskip
\noindent En dimension 2 ou si $\mathcal{G}$ est complexe, ce résultat peut se reformuler comme 
un théorème de densité des cocycles réductibles au voisinage des cocycles constants:

\begin{thm}\label{th3}Soit $\mathcal{G}=gl(n,\mathbb{C}),u(n), gl(2,\mathbb{R}), sl(2,\mathbb{R})$ ou $o(2)$. Soit $0<r\leq 1$ et $A\in \mathcal{G}, F\in C^{G,2}_r(\mathbb{T}^d,\mathcal{G})$. Il existe 
$\epsilon_0$ ne dépendant que de $r,n,d,\kappa,\tau,A$ tel que si 

$$| |F| |_r\leq \epsilon_0$$ 

\noindent alors pour tout 
$\epsilon>0$ il existe $r_\epsilon>0,H\in C^\omega_{r_\epsilon}(\mathbb{T}^d,\mathcal{G})$ réductible dans $C^{G,2}_{r_\epsilon}(\mathbb{T}^d,\mathcal{G})$ et tel que 

$$| | A+F-H | |_{r_\epsilon}\leq \epsilon$$ 
\end{thm}


\subsection{Résultats déjà connus}

\noindent Un résultat comparable pour les cocycles lisses à valeurs dans les groupes compacts avait été obtenu par R. Krikorian dans \cite{Kr} (th.5.1.1). Dans le cas d'un cocycle au-dessus d'une rotation du cercle, le contrôle de l'analyticité est bien meilleur (voir par exemple \cite{AK06}) car il est alors possible d'utiliser des méthodes globales. Nous considérons ici le cas d'un tore de dimension quelconque. La méthode KAM que nous allons utiliser ici avait déjà produit des résultats de réductibilité en mesure totale pour des cocycles à valeurs dans $SL(2,\mathbb{R})$ (\cite{El92}, \cite{SH06}).

\bigskip 
\noindent Un résultat analogue au théorème \ref{th1} dans le cadre des cocycles analytiques à valeurs dans $GL(n,\mathbb{R})$ avait déjà été démontré 
dans \cite{E1} par L.H.Eliasson; on se reportera également à sa généralisation dans \cite{C2}. Notons bien, cependant, la distinction suivante: le théorème 
\ref{th1} est un résultat de presque réductibilité faible, c'est-à-dire que le paramètre peut tendre vers 0; il en est de même pour le théorème d'Eliasson dans le cas des fonctions analytiques (\cite{E1}) où le rayon d'analyticité peut tendre vers 0. Dans le cas des fonctions analytiques, on peut obtenir un résultat de presque réductibilité forte (\cite{C2}), c'est-à-dire où le rayon d'analyticité reste strictement positif à la limite. Mais la question de la presque réductibilité forte dans le cas des fonctions Gevrey reste ouverte.







\bigskip 
\noindent Notons que, comme dans \cite{E1} et \cite{C2}, la perte de périodicité dans le théorème \ref{th1} est inévitable dans un groupe 
réel de dimension plus grande que 2. La notion de "bonnes propriétés de périodicité" a pour but de s'assurer qu'un 
seul doublement de période suffit. Le cadre symplectique introduit de nouvelles contraintes pour éliminer les 
résonances, par rapport au cadre réel, mais ces contraintes n'ont pas de conséquences sur la construction de la fonction de renormalisation, ce qui explique 
que l'on n'ait pas plus de perte de périodicité que dans le cadre $GL(n,\mathbb{R})$. Ainsi, comme dans \cite{C}, un seul doublement de période suffit dans le cas d'un groupe symplectique réel.

\subsection{Plan de l'article}

La preuve des théorèmes \ref{th1} et \ref{th3} reprend à peu près le même schéma de démonstration que dans \cite{E1}; il s'agit d'une preuve par itération 
de type KAM. En voici les principales étapes:

\bigskip

\begin{itemize}

\item Construction d'une renormalisation $\Phi$ d'ordre $N$ (proposition \ref{renormG})
 pour $N\in \mathbb{N}\setminus\{0\}$ bien choisi. 

\bigskip
\noindent Notons qu'en dimension 2, $\Phi$ est telle que pour toute fonction $H$ continue sur 
$\mathbb{T}^d$, $\Phi H\Phi^{-1}$ est continue sur $\mathbb{T}^d$.

\bigskip

\item Résolution de l'équation homologique (proposition \ref{homolG}): si 
$\tilde{A}$ a un spectre vérifiant certaines conditions diophantiennes 
et que $\tilde{F}$ est 
une fonction qui a certaines bonnes propriétés de périodicité relatives à $\tilde{A}$, alors il existe une solution 
$\tilde{X}$ de l'équation 

$$\partial_\omega \tilde{X}=[\tilde{A},\tilde{X}]+\tilde{F}^{{N}};\ \hat{\tilde{X}}(0)=0$$

\noindent qui a les mêmes propriétés de périodicité que $\tilde{F}$; elle prend ses valeurs dans la même algèbre de Lie que $\tilde{F}$. 
De plus, elle vérifie une bonne estimation 
quitte à perdre un peu de régularité.

\bigskip

\item Lemme inductif (proposition \ref{iter2G}):
Si $\tilde{F}\in C^{G,2}_{r}(2\mathbb{T}^d,\mathcal{G})$ a certaines propriétés de périodicité (relatives à $
\tilde{A}$), 
si 

$$\partial_\omega \Psi=\bar{A}\Psi-\Psi \tilde{A} $$

\noindent et $\bar{F}=\Psi \tilde{F} \Psi^{-1}$, alors il existe 
$Z\in C^{G,2}_{r'}(2\mathbb{T}^d,G)$ telle que 

\begin{equation}\partial_\omega Z =(\bar{A} +\bar{F} )Z -Z (\bar{A}' 
+\bar{F}' )
\end{equation}

\noindent où $\bar{A}'$ est réductible, $\bar{F}'$ est beaucoup plus petit que $\bar{F}$, 
$Z$ est proche de l'identité et $\Psi'^{-1}\bar{F}'\Psi'$ a des propriétés de périodicité relatives à $A'$ 
analogues à celles de $\tilde{F}$.

\bigskip
\noindent L'estimation de $\bar{F}'$ dépend de $\tilde{F}-\tilde{F}^N$, de la fonction de renormalisation $\Phi$, 
et de la solution $\tilde{X}$ de l'équation homologique. 

\bigskip

\item Itération du lemme inductif (théorème \ref{PRG}):
On itèrera le lemme \ref{iter2G} grâce à un lemme numérique (lemme \ref{eps'1G}), 
pour réduire arbitrairement la perturbation.

\bigskip

\rem La fonction de renormalisation se construit selon le même principe que dans \cite{E1}.

\end{itemize}

\subsection{Notations}

\noindent 
On notera $\langle .,. \rangle$ le produit scalaire euclidien complexe, avec la convention qu'il est antilinéaire en la deuxième variable.
Pour tout opérateur linéaire $M$, on notera $M^*$ son adjoint, égal à la transposée de $M$ dans le cas où $M$ est réel. On notera aussi $M_{\mathcal{N}}$ la partie nilpotente de 
$M$, c'est-à-dire:
si $M=PAP^{-1}$ avec $A$ en forme normale de Jordan, et si $A_D$ représente la partie diagonale de $A$, alors $M_{\mathcal{N}}=P(A-A_D)P^{-1}$.
Pour simplifier l'écriture, si $A: 2\mathbb{T}^d\rightarrow GL(n,\mathbb{C})$, on notera 
$A^{-1}$ la fonction ${\theta\mapsto A(\theta)^{-1}}$. Pour tout $m=(m_1,\dots ,m_d)\in \frac{1}{2}\mathbb{Z}^d$, on notera $\mid m\mid = \mid m_1\mid +\dots +\mid m_d\mid$. 
On désignera par $J$ la matrice $
J=\left(\begin{array}{cc}
0 & -Id\\
Id & 0\\
\end{array}\right)
$.



\section{Bonnes propriétés de périodicité}

\noindent Commençons par introduire quelques définitions. La notion de trivialité par rapport à une décomposition permettra 
d'expliciter plus facilement les fonctions de renormalisation; les bonnes propriétés de périodicité ont déjà été 
introduites dans \cite{E1} et permettent de s'assurer, dans le cas réel, qu'un seul doublement de
 période est nécessaire à l'itération du lemme inductif. 

\subsection{Décompositions invariantes}

\noindent \Def $\mathcal{L}=\{L_1,\dots ,L_R\}$ est une \textit{décomposition de $\mathbb{C}^n$} si $\mathbb{C}^n=\bigoplus_j L_j$. 

\bigskip
\noindent \Def Soient $\mathcal{L},\mathcal{L}'$ des décompositions de $\mathbb{C}^n$. On dit que $\mathcal{L}$ est 
\textit{plus fine que} 
$\mathcal{L}'$ si pour tout $L\in\mathcal{L}$, il existe $L'\in\mathcal{L}'$ tel que $L\subset L'$; 
on dit que $\mathcal{L}$ est \textit{strictement plus fine que} 
$\mathcal{L}'$ si $\mathcal{L}$ est plus fine que 
$\mathcal{L}'$ et $\mathcal{L}\neq\mathcal{L}'$.

\bigskip
\noindent \Def Soit $A\in gl(n,\mathbb{C})$; on dit que $\mathcal{L}=\{L_1,\dots ,L_s\}$ est une \textit{$A$-décomposition}, ou \textit{décomposition 
$A$-invariante}, 
si 
c'est une décomposition de $\mathbb{C}^n$ et que pour tout $i$, $AL_i\subset L_i$. Les $L_i$ sont les \textit{sous-espaces} de 
$\mathcal{L}$. 


\bigskip
\noindent \rem Une $A$-décomposition est toujours moins fine que la décomposition en sous-espaces propres généralisés. Ainsi, si deux matrices $A$ et $A'$ ont la même 
décomposition en sous-espaces propres généralisés, alors une $A$-décomposition est une $A'$-décomposition.

\bigskip
\noindent \textbf{Notation:} Soit $\mathcal{L}$ une $A$-décomposition. Pour tout $L\in\mathcal{L}$, on note 
$\sigma(A_{|L})$ le spectre de $A$ restreint au sous-espace $L$.

\bigskip
\noindent \Def Soit $\kappa'\geq 0$. On notera $\mathcal{L}_{A,\kappa'}$ l'unique $A$-décomposition $\mathcal{L}$ telle que pour tous $L\neq L'\in \mathcal{L}$, $\alpha\in\sigma(A_{|L})$ et 
$\beta\in\sigma(A_{|L'})$ 
$\Rightarrow |\alpha-\beta|>\kappa'$ et qu'aucune $A$-décomposition strictement plus fine que $\mathcal{L}$ n'a cette propriété.

\bigskip
\noindent \Def Soit $A\in gl(n,\mathbb{C})$. 
On note $\mathcal{L}_A$ la décomposition de $\mathbb{C}^n$ qui est l'ensemble des sous-espaces propres généralisés de $A$.

\bigskip
\noindent \rem On fera attention au fait que $\mathcal{L}_A$ n'est pas forcément identique à $\mathcal{L}_{A,0}$. En général $\mathcal{L}_A$ est plus fine.

\bigskip
\noindent \Def Soit $\mathcal{L}$ une décomposition de $\mathbb{C}^n$. Pour tout $u\in\mathbb{C}^n$, il existe une unique 
décomposition $u=\sum_{L\in\mathcal{L}}u_L$ avec $u_L\in L$ pour tout $L\in\mathcal{L}$. 
Pour tout $L\in\mathcal{L}$, on appelle \textit{projection 
sur $L$ relative à $\mathcal{L}$}, et on note $P^\mathcal{L}_L$, l'application définie par $P^\mathcal{L}_Lu=u_L$.

\bigskip
\noindent 
\rem Soit $A\in gl(n,\mathbb{C})$ et $\kappa'>0$. Si $\mathcal{L}$ est une $A$-décomposition moins fine que $\mathcal{L}_{A,\kappa'}$, alors on a le lemme suivant, tiré de \cite{E1}, appendice, lemme A\footnote{Le lemme A de \cite{E1} donne en fait une estimation en fonction de $\mid \mid A\mid \mid$, mais il apparaît clairement dans sa démonstration que l'estimation ne dépend en fait que de $A_{\mathcal{N}}$.}:

\begin{lem}\label{c0}
Il existe une constante $C_0\geq 1$ ne dépendant que de $n$ telle que tout sous-espace $L\in \mathcal{L}$ vérifie

\begin{equation}\mid \mid P^{\mathcal{L}}_L\mid \mid \leq C_0\left(\frac{1+\mid \mid A_{\mathcal{N}}\mid \mid }{\kappa'}\right)^{n(n+1)}
\end{equation}

\end{lem}

\noindent Par la suite, on notera toujours $C_0$ cette constante fixée dans le lemme \ref{c0}.

\bigskip
\Def Une \textit{$(A,\kappa',\gamma)$-décomposition} est une $A$-décomposition $\mathcal{L}$ telle que
%
pour tout $L\in \mathcal{L}$, la projection sur $L$ relative à $\mathcal{L}$ 
vérifie l'estimation

\begin{equation}\mid \mid P^{\mathcal{L}}_L\mid \mid \leq C_0\left(\frac{1+\mid \mid A_{\mathcal{N}}\mid \mid }{\kappa'}\right)^{\gamma}
\end{equation}


\bigskip
\noindent \rem 
 Pour $A\in gl(n,\mathbb{C})$, on a toujours $A=\sum_{L,L'\in\mathcal{L}}P^\mathcal{L}_L AP^\mathcal{L}_{L'}$. 
En particulier, si $\mathcal{L}$ est une $A$-décomposition, alors $A=\sum_{L\in\mathcal{L}}P^\mathcal{L}_L A P^\mathcal{L}_L$.

\bigskip
\noindent \Defs Soit $\mathcal{L}$ une décomposition. 
\begin{itemize}
\item On dit que $\mathcal{L}$ est une \textit{décomposition réelle} si 
pour tout 
$L\in\mathcal{L}$, $\bar{L}\in\mathcal{L}$;

\item On dit que $\mathcal{L}$ est une \textit{décomposition symplectique} si c'est une décomposition de $\mathbb{C}^{n}$ avec $n$ pair et que pour tout $L\in \mathcal{L}$, il existe un unique $L'
\in \mathcal{L}$ tel que $\langle L,JL'\rangle \neq 0$. 
\item On dit que $\mathcal{L}$ est une \textit{décomposition unitaire} si pour tout $L\neq L'\in \mathcal{L}$, 
$\langle L,L'\rangle =0$.

\end{itemize}

\bigskip
\noindent \rem 
\begin{itemize}
\item Si $A$ est une matrice réelle, alors pour tout $\kappa'\geq 0$,  $\mathcal{L}_{A,\kappa'}$ est une décomposition réelle. 
\item Pour tout $L$, il existe au moins un $L'$ tel que $\langle L,JL'\rangle \neq 0$. Cela vient de la non-dégénérescence de la forme symplectique $\langle .,J.\rangle$.
\item Si $A\in sp(n,\mathbb{R})$, alors toute $A$-décomposition $\mathcal{L}$ moins fine que $\mathcal{L}_{A,0}$ est une décomposition symplectique réelle. En effet, soient $L,L'\in \mathcal{L}$ tels que $\langle L,JL'\rangle \neq 0$; soient $v\in L,v'\in L'$ des vecteurs propres de $A$ tels que $\langle v,Jv'\rangle \neq 0$ 
et $\lambda,\lambda'$ leurs valeurs propres associées. Alors 

\begin{equation}\lambda \langle v,Jv'\rangle= \langle Av,Jv'\rangle =\langle v,A^*Jv'\rangle 
=-\langle v,JAv'\rangle =-\bar{\lambda}'\langle v,Jv'\rangle 
\end{equation}

\noindent et comme $\langle v,Jv'\rangle \neq 0$, alors $\lambda=-\bar{\lambda}'$.

\item Si $A\in U(n)$, alors toute décomposition moins fine que $\mathcal{L}_{A,0}$ est une décomposition unitaire. 

\end{itemize}

\subsection{Trivialité et bonnes propriétés de périodicité par rapport à une décomposition}

\noindent \Def Soit $\mathcal{L}$ une décomposition de $\mathbb{C}^n$. 
On dit qu'une fonction $\Psi$ 
est \textit{triviale par rapport à $\mathcal{L}$} s'il existe 
${\{m_L,\ L\in \mathcal{L}\}\subset \frac{1}{2}\mathbb{Z}^d}$, tels que pour tout $\theta\in 2\mathbb{T}^d$,

\begin{equation}
\Psi(\theta) =\sum_{L\in\mathcal{L}}e^{2i\pi\langle m_L,\theta\rangle}P_L^\mathcal{L}\end{equation}

\bigskip
\noindent \Def On dit que la fonction $\Psi$ est \textit{triviale} s'il existe une décomposition $\mathcal{L}$ telle que $\Psi$ est triviale par rapport à $\mathcal{L}$. 

\bigskip
\noindent \rem 
\begin{itemize}

\item Si $\Psi$ est triviale 
par rapport à $\mathcal{L}$ et que $\mathcal{L}'$ est plus fine que $\mathcal{L}$, 
alors elle est triviale par rapport à $\mathcal{L}'$. 
\item Si $\Phi,\Psi:2\mathbb{T}^d\rightarrow GL(n,\mathbb{C})$ sont triviales par rapport à 
$\mathcal{L}$, alors 
le produit $\Phi\Psi$ est trivial par rapport à $\mathcal{L}$.
\item Si $\Phi$ est triviale par rapport à une $A$-décomposition $\mathcal{L}$, alors pour tout $\theta\in
2\mathbb{T}^d$, $[A,\Phi(\theta)]=0$.

\end{itemize}

\begin{lem}\label{réel} Soit $\mathcal{L}$ une décomposition réelle de $\mathbb{C}^n$, ${\{m_L,\ L\in \mathcal{L}\}\subset \frac{1}{2}\mathbb{Z}^d}$ 
et $\Psi$ définie par

\begin{equation}
\label{psitriv}
\Psi(\theta) =\sum_{L\in\mathcal{L}}e^{2i\pi\langle m_L,\theta\rangle}P_L^\mathcal{L}\end{equation}

\noindent Alors $\Psi$ est réelle si et seulement si pour tout $L$, $m_L=-m_{\bar{L}}$. De plus, si $\Psi$ est réelle, alors $\Psi$ est à valeurs dans $SL(n,\mathbb{R})$.

\end{lem}

\noindent Pour la démonstration, se reporter à \cite{C2}, lemme 1.2.

\bigskip
\noindent \rem Toute fonction triviale par rapport à une décomposition unitaire est unitaire. En effet, soit $\mathcal{L}$ une décomposition unitaire, $\Phi$ triviale par 
rapport à $\mathcal{L}$ et $L,L'\in \mathcal{L}$. 
Alors pour tout $u\in \mathcal{L},v\in \mathcal{L}'$,

\begin{equation}\langle \Phi(\theta) u,\Phi(\theta) v\rangle 
= \langle e^{2i\pi \langle m_L,\theta \rangle} u,e^{2i\pi \langle m_{L'},\theta \rangle} v\rangle 
=\langle  u, v\rangle 
\end{equation}

\begin{lem}\label{phisp}Soit $\mathcal{L}$ une décomposition symplectique réelle et $\{m_L,L\in\mathcal{L}\}$ une famille d'éléments de $\frac{1}{2}\mathbb{Z}^d$. 
Soit $\Psi=\sum_{L\in \mathcal{L}}e^{2i\pi \langle m_L,.\rangle}P^{\mathcal{L}}_L$.
Alors $\Psi $ est à valeurs dans $Sp(n,\mathbb{R})$ si et seulement si 
\begin{itemize}
\item pour tout $L$, $m_L=-m_{\bar{L}}$
\item et si $\langle L,JL'\rangle \neq 0$, 
alors $m_L=m_{L'}$. 
\end{itemize}

\end{lem}

\noindent \dem Se référer à \cite{C2}, lemme 1.3.

\bigskip
\noindent Définissons les propriétés de périodicité. 

\bigskip
\noindent \Def Soit $\mathcal{L}$ une décomposition de $\mathbb{C}^n$. 
On dit que $F\in C^0(2\mathbb{T}^d,
gl(n,\mathbb{R}))$ \textit{a de bonnes propriétés de périodicité par rapport à $\mathcal{L}$} s'il existe $\Phi$ triviale par rapport à $\mathcal{L}$ telle que 
$\Phi^{-1}F\Phi$ soit continue sur $\mathbb{T}^d$.

\noindent Pour préciser la famille $(m_L)$ qui définit $\Phi$, on dira que $F$ \textit{a de bonnes propriétés de périodicité par rapport à 
$\mathcal{L}$ 
et $(m_L)$}.

\bigskip
\noindent \rem 
\begin{itemize}
\item Si $F\in C^0(2\mathbb{T}^d,gl(n,\mathbb{R}))$ a de bonnes propriétés de périodicité par rapport à une décomposition $\mathcal{L}$ et que $\Phi$ est 
triviale par rapport à $\mathcal{L}$, 
alors $\Phi F \Phi^{-1}$ a de bonnes propriétés de périodicité par rapport 
à $\mathcal{L}$.

\item 
Si $\mathcal{L}'$ est une décomposition de $\mathbb{C}^n$ qui est plus fine que $\mathcal{L}$ et que $F$ a de bonnes propriétés de périodicité par rapport à $\mathcal{L}$, alors 
$F$ a de bonnes propriétés de périodicité par rapport à $\mathcal{L}'$. 

\item Soit 
$\mathcal{L}$ une décomposition de $\mathbb{C}^n$ et $(m_L)_{L\in \mathcal{L}}$ une famille d'éléments de $\frac{1}{2}\mathbb{Z}^d$. 
Si $F_1,F_2\in C^0(2\mathbb{T}^d, gl(n,\mathbb{R}))$ ont de bonnes propriétés de 
périodicité par rapport à $\mathcal{L}$ et $(m_L)$, alors le produit $F_1F_2$ a de bonnes propriétés de 
périodicité par rapport à $\mathcal{L}$ et $(m_L)$.

\end{itemize}

\section{Généralités sur les fonctions de classe Gevrey deux}

\rem \begin{itemize}
\item Pour tous $0<r'<r$, on a l'inclusion $C^{G,2}_r( 2\mathbb{T}^d,\mathcal{G})\subset C^{G,2}_{r'}( 2\mathbb{T}^d,\mathcal{G})$ et $\mid \mid f \mid \mid_{r'}
\leq \mid \mid f \mid \mid_r$. 

\item Pour $f,g\in C^{G,2}_r(2\mathbb{T}^d,\mathcal{G})$, on a $||fg||_r\leq ||f||_r||g||_r$ (voir par exemple 
\cite{MS}, appendice).
\end{itemize}

\begin{lem}\label{gevreyexp} Pour tout $m\in\mathbb{Z}^d$ et tout $r'>0$, la fonction $\theta\mapsto e^{2i\pi \langle m,\theta\rangle}$ vérifie 

\begin{equation}\mid \mid e^{2i\pi \langle m,.\rangle}\mid \mid _{r'}\leq  e^{2\pi|m|r'^2}
\end{equation}

 \end{lem}

\dem
Pour tout $\alpha\in \mathbb{N}^d$ et tout $\theta\in \mathbb{T}^d$,

\begin{equation}\begin{split}\frac{r'^{2\mid \alpha\mid }}{(\alpha!)^2} \mid \partial^\alpha (e^{2i \pi \langle m,\theta \rangle})\mid 
&\leq \frac{r'^{2\mid \alpha\mid }}{(\alpha!)^2}  \prod_j \mid 2\pi m_j\mid ^{\alpha_j} \\
&\leq \prod_j \frac{( r'^2\mid 2\pi m_j\mid )^{ \alpha_j }}{(\alpha_j!)^2}
\leq \prod_j e^{| 2\pi m_j|r'^2}=e^{2\pi |m|r'^2}\ \Box
\end{split}\end{equation}
%







\noindent \rem Ceci implique que les fonctions analytiques sur un $r$-voisinage du tore ou du double tore sont Gevrey deux de paramètre $r$;

\section{Fonction de renormalisation}


\noindent Rappelons un lemme de \cite{C2} permettant de construire une fonction de renormalisation à l'ordre $N$:

\begin{prop}\label{diophanG}Soit $A\in \mathcal{G}$ et $N\in\mathbb{N}$. 
Soit

\begin{equation}\kappa''= \frac{\kappa}{n(8{N})^{\tau}}\end{equation}

\noindent 
Il existe 
une fonction $\Phi$ triviale par rapport à $\mathcal{L}_{A,\kappa''}$ et à valeurs dans $G$ telle que

\begin{enumerate}

\item pour tout $r'\geq 0$,

\begin{equation}\label{normephi}|\Phi|_{r'}\leq nC_0\left(\frac{1+||A_{\mathcal{N}}||}{\kappa''}\right)^{n(n+1)}e^{4\pi {N}r'},\ |\Phi^{-1}|_{r'}\leq 
nC_0\left(\frac{1+||A_{\mathcal{N}}||}{\kappa''}\right)^{n(n+1)}e^{4\pi {N}r'}
\end{equation}

\item Si $\tilde{A}$ est définie par 

\begin{equation}\forall \theta\in 2\mathbb{T}^d,\ \partial_\omega \Phi(\theta)=A\Phi(\theta)-\Phi(\theta)\tilde{A}
\end{equation}

\noindent alors 

\begin{equation}\label{normeA} ||\tilde{A}-A||\leq 4\pi {N}\end{equation}

\noindent et $\tilde{A}$ a un spectre $DC^{N}_\omega (\kappa'',\tau)$.

\item Si $\mathcal{G}=gl(n,\mathbb{C})$ ou $u(n)$, $\Phi$ est définie sur $\mathbb{T}^d$.

\end{enumerate}

\end{prop}

\Def On appellera la fonction $\Phi$ construite dans la proposition \ref{diophanG} une \textit{renormalisation de $A$ d'ordre $N$}. 

\bigskip
\noindent 
Cette fonction vérifie une bonne estimation en norme Gevrey, comme le montre la proposition suivante:

\begin{lem}\label{renormG} Soit $N\geq 2,A\in gl(n,\mathbb{C})$ et $\Phi$ une renormalisation de $A$ d'ordre ${N}$. 
Alors $\Phi$ vérifie, pour 
tout $r'$, 
l'estimation en norme Gevrey

\begin{equation}\label{normephiG}||\Phi||_{r'}\leq nC.C_0 \left(\frac{1+||A_{\mathcal{N}}||}{\kappa''}\right)^{n(n+1)} 
e^{4\pi r'^2{{N}}}
\end{equation}

\noindent où $C$ ne dépend que de $d$, et de même pour $\Phi^{-1}$. De plus, si $\mathcal{G}=o(n)$ ou $u(n)$, alors

\begin{equation}\label{normephiG'}||\Phi||_{r'}\leq nC
e^{4\pi r'^2{{N}}}
\end{equation}

\noindent et de même pour $\Phi^{-1}$.
\end{lem}

\noindent \dem Pour tout $m\in \mathbb{Z}^d$ et tout $r'>0$, d'après le lemme \ref{gevreyexp},

\begin{equation}\begin{split}||e^{2i\pi\langle m,.
\rangle}||_{r'}
\leq Ce^{2\pi r'^2 |m|}
\end{split}\end{equation}

\noindent où $C$ ne dépend que de $d$. Donc

\begin{equation}\begin{split}
||\Phi||_{r'}&\leq \sum_{L\in \mathcal{L}_{A,\kappa''}} \mid \mid P^{\mathcal{L}_{A,\kappa''}}_L \mid \mid \ \mid\mid e^{2i\pi \langle m_L,.\rangle}\mid\mid_{r'}\\
&\leq C \sum_{L\in \mathcal{L}_{A,\kappa''}} \mid \mid P^{\mathcal{L}_{A,\kappa''}}_L \mid \mid e^{2\pi r'^2 |m_L|}\\
&\leq C \sum_{L\in \mathcal{L}_{A,\kappa''}} \mid \mid P^{\mathcal{L}_{A,\kappa''}}_L \mid \mid e^{4\pi r'^2 {N}}
\end{split}\end{equation}

\noindent or d'après le lemme \ref{c0}, 

\begin{equation}\mid \mid P^{\mathcal{L}_{A,\kappa''}}_L \mid \mid \leq 
C_0\left(\frac{1+\mid \mid A_{\mathcal{N}}\mid \mid }{\kappa''}\right)^{n(n+1)}
\end{equation}

\noindent d'où \eqref{normephiG}. 
Si $\mathcal{G}$ est $o(n) $ ou $u(n)$, alors $\mathcal{L}_{A,\kappa''}$ est une décomposition unitaire et donc $P^{\mathcal{L}_{A,\kappa''}}_L$ est de norme 1, 
d'où \eqref{normephiG'}. 
$\Box$

\section{Equation homologique}


\begin{lem}\label{gevrey}Soit $0<r'<r\leq 1$, $f\in C^{G,2}_r(2\mathbb{T}^d,\mathcal{G})$ et 
$g\in C^{G,2}_{r'}(2\mathbb{T}^d,\mathcal{G})$. Soient ${C>0,D\geq 0}$. Supposons que pour tout $m\in \frac{1}{2}\mathbb{Z}^d$, 

\begin{equation}||\hat{g}(m)||\leq C |m|^{D}||\hat{f}(m)||\end{equation}

\noindent 
Alors

\begin{equation}||g||_{r'}\leq C'C ||f||_r\left(\frac{1}{r-r'}\right)^{3(D+2)}
\end{equation}

\noindent où $C'$ ne dépend que de $d,D$.
\end{lem}

\dem On a pour tout $\theta\in N\mathbb{T}^d$ et tout $\alpha\in \mathbb{N}^d$,

\begin{equation}\begin{split}||\partial^\alpha g(\theta)||&\leq \sum_{m\in \frac{1}{2}\mathbb{Z}^d} \mid\mid \hat{g}(m)\mid \mid \ \mid  \partial^\alpha e^{2i\pi \langle m,\theta\rangle}    \mid\\
&\leq \sum_{m\in \frac{1}{2}\mathbb{Z}^d}  \mid\mid \hat{g}(m) \mid \mid \ \prod_j \mid 2\pi m_j\mid^{\alpha_j}\\
\end{split}\end{equation}

\noindent donc par hypothèse, 

\begin{equation}\begin{split}||\partial^\alpha g(\theta)||
&\leq C\sum_{m\in \frac{1}{2}\mathbb{Z}^d}  |m|^{D}||\hat{f}(m)||    \prod_j \mid 2\pi m_j\mid^{\alpha_j}\\
&\leq C'C \sum_{m\neq 0} \frac{1}{\mid 2\pi m\mid ^{2d}}
\prod_i
\mid 2\pi m_i\mid ^{\alpha_i+(D+2)e_i} ||\hat{f}(m)||\\
\end{split}\end{equation}

\noindent et donc 
en notant $\bar{1}=(1,\dots,1)$,

\begin{equation}\begin{split}||\partial^\alpha g(\theta)||
&\leq C'C \sum_{m\neq 0} 
\frac{1}{\mid 2\pi m\mid ^{2d}}
 ||\widehat{\partial^{\alpha+(D+2)\bar{1}}f}(m)||\\
&\leq C' C\sup_\theta ||\partial^{\alpha+(D+2)\bar{1}}f(\theta)||
\end{split}\end{equation}



\noindent où $C'$ ne dépend que de $d,D$, donc 

\begin{equation}\begin{split}||g||_{r'}&=\sup_{\alpha\in \mathbb{N}^d} (r')^{2\mid \alpha\mid}\frac{1}{(\alpha!)^2}\sup_\theta ||\partial^\alpha g(\theta)||\\
&\leq C'C\sup_{\alpha} (r')^{2\mid \alpha\mid }\frac{1}{(\alpha!)^2} \sup_{\theta} ||\partial^{\alpha+(D+2)\bar{1}}f(\theta)||\\
& \leq C' C\sup_{\alpha} \frac{r^{2\mid \alpha+(D+2)\bar{1}\mid }}{((\alpha+(D+2)\bar{1})!)^2} 
\sup_\theta  ||\partial^{\alpha+(D+2)\bar{1}} f(\theta)||
\frac{r'^{2\mid \alpha\mid}}{r^{2\mid \alpha+(D+2)\bar{1}\mid }}\left(\frac{(\alpha+(D+2)\bar{1})!}{\alpha!}\right)^2\\
& \leq C'C  \mid \mid f\mid \mid _r 
\sup_{\alpha}\frac{r'^{2\mid \alpha\mid}}{r^{2(\mid \alpha\mid +(D+2)d) }}\prod_i\left(\frac{(\alpha_i+D+2)!}{\alpha_i!}\right)^2\\
& \leq C'C  \mid \mid f\mid \mid _r 
\sup_{\alpha}\frac{r'^{2\mid \alpha\mid}}{r^{2(\mid \alpha\mid +(D+2)d) }}\left(\frac{(\mid \alpha\mid +D+2)!}{\mid \alpha\mid !}\right)^{2d}\\
& \leq C' C \mid \mid f\mid \mid _r 
\sup_{\alpha}\frac{r'^{2\mid \alpha\mid}}{r^{2(\mid \alpha\mid +(D+2)d) }}\left( \mid \alpha\mid +D+2\right)^{2(D+2)d}\\
\end{split}\end{equation}

\noindent Or la fonction 

$$\phi:[0,+\infty[\rightarrow [0,+\infty[
,\ t\mapsto\left( \frac{r'}{r}\right)^t t^{2(D+2)d}$$ 

\noindent atteint son maximum en $t=\frac{2(D+2)d}{\ln \frac{r}{r'}}$ où elle vaut 
$e^{-2(D+2)d}\left(\frac{2(D+2)d}{\ln \frac{r}{r'}}\right)^{2(D+2)d}$. D'où

\begin{equation}\begin{split}||g||_{r'}&\leq C'C ||f||_re^{-2(D+2)d}\left(\frac{2(D+2)d}{{r'}\ln \frac{r}{r'}}\right)^{2(D+2)d}\\
&\leq C'C ||f||_re^{-2(D+2)d}\frac{(2(D+2)d)^{2(D+2)d}}{(r-r')^{3(D+2)d}}\ \Box
\end{split}\end{equation}

\begin{prop}\label{homolG} 
Soient 
\begin{itemize}
\item ${N}\in\mathbb{N}$,
\item $\kappa'\in ]0, \kappa]$, 
\item $\gamma \geq n(n+1)$,
\item
 $0<r'<r$. 
 \end{itemize}
Soit $\tilde{A}\in \mathcal{G}$ avec un spectre $DC^{{N}}_\omega(\kappa',\tau)$.
Soit $\tilde{F}\in C^\omega_r(2\mathbb{T}^d,  \mathcal{G})$ avec de bonnes propriétés de périodicité par rapport à une $(\tilde{A},\kappa',\gamma)$-décomposition $\mathcal{L}$. 
Alors il existe une solution $\tilde{X}\in C^\omega_{r'}(2\mathbb{T}^d,  \mathcal{G})$ de l'équation 

\begin{equation}\label{hom}\forall \theta\in 2\mathbb{T}^d, \ \partial_\omega \tilde{X}(\theta)=[\tilde{A},\tilde{X}(\theta)]+\tilde{F}^{{N}}(\theta)
-\hat{\tilde{F}}(0);\ \hat{\tilde{X}}(0)=0
\end{equation}

\noindent telle que 
\begin{itemize}
\item si $\tilde{F}$ a de bonnes propriétés de périodicité par rapport à 
$\mathcal{L}$ et $(m_L)$, 
alors $\tilde{X}$ a de bonnes propriétés de périodicité par rapport à $\mathcal{L}$ et $(m_L)$; en particulier, si $\tilde{F}$ est définie sur $\mathbb{T}^d$, alors $\tilde{X}$ l'est également,

\item  Il existe $C',D'$ ne dépendant que de 
$n,d,\tau$ tels que,

\begin{equation}\label{X2G}\mid \mid \tilde{X}\mid \mid _{r'}
\leq C'\left(\frac{1+||\tilde{A}_{\mathcal{N}}||}{(r-r')\kappa'}\right)^{D'\gamma}
\mid \mid \tilde{F}\mid \mid _r\end{equation}

\end{itemize}

\bigskip
De plus, la troncation de 
$\tilde{X}$ à l'ordre $N$ est unique. 

\end{prop}

%


%

\dem $\bullet$ L'existence de $\tilde{X}$, son unicité jusqu'à l'ordre $N$, le fait qu'elle prend ses valeurs dans $\mathcal{G}$ et ses bonnes propriétés de périodicité par rapport à $\mathcal{L}$ se démontrent comme dans \cite{C2}, proposition 3.2.

\bigskip
\noindent $\bullet$ Pour obtenir l'estimation \eqref{X2G}, on démontre d'abord que pour tout $m\in \frac{1}{2}\mathbb{Z}^d$ et tous $L,L'\in \mathcal{L}'$,

\begin{equation}\label{estimationreG}||P^{\mathcal{L}}_L\hat{\tilde{X}}(m)
P^{\mathcal{L}}_{L'}||
\leq C'\frac{(1+||\tilde{A}_{\mathcal{N}}||)^{n^2-1}|m|^{(n^2-1)\tau}}{\kappa'^{(n^2-1)}}
||P^{\mathcal{L}}_L\hat{\tilde{F}}(m)
P^{\mathcal{L}}_{L'}||  (||P^{\mathcal{L}}_L|| \ 
||P^{\mathcal{L}}_{L'}||)^{n^2-1}
\end{equation}

\noindent Pour cela, nous allons faire un raisonnement analogue à celui de \cite{E1}, lemme 2. Soit $\mathcal{A}_{L,L'}$ l'opérateur linéaire allant de $gl(n,\mathbb{C})$ dans lui-même 
tel que pour tout $M\in gl(n,\mathbb{C})$,

\begin{equation}\mathcal{A}_{L,L'}M=
\tilde{A} P^{\mathcal{L}}_L
M-M P^{\mathcal{L}}_{L'}
 \tilde{A} 
\end{equation}

\noindent En décomposant \eqref{hom} en blocs, on obtient pour tous $L,L'\in \mathcal{L}$

\begin{equation}\begin{split}&\partial_\omega (P^{\mathcal{L}}_L \tilde{X}P^{\mathcal{L}}_{L'})
=\mathcal{A}_{L,L'} P^{\mathcal{L}}_L\tilde{X}P^{\mathcal{L}}_{L'}
+P^{\mathcal{L}}_L(\tilde{F}^{{N}}-\hat{\tilde{F}}(0))P^{\mathcal{L}}_{L'}
\end{split}\end{equation}

\noindent Donc pour tout $m\in \frac{1}{2}\mathbb{Z}^d$ tel que $0<\mid m\mid \leq N$,

\begin{equation}\label{fourierAR}2i\pi\langle m,\omega\rangle (P^{\mathcal{L}}_L \hat{\tilde{X}}(m)
P^{\mathcal{L}}_{L'})
=\mathcal{A}_{L,L'}(P^{\mathcal{L}}_L \hat{\tilde{X}}(m)
P^{\mathcal{L}}_{L'})
+P^{\mathcal{L}}_L\hat{\tilde{F}}(m)P^{\mathcal{L}}_{L'}
\end{equation}

donc

\begin{equation}\label{invR} (P^{\mathcal{L}}_L \hat{\tilde{X}}(m)
P^{\mathcal{L}}_{L'})
=(2i\pi\langle m,\omega\rangle-\mathcal{A}_{L,L'})^{-1}P^{\mathcal{L}}_L\hat{\tilde{F}}(m)P^{\mathcal{L}}_{L'}
\end{equation}

Représentons $\mathcal{A}_{L,L'}$ comme une matrice de dimension $n^2$. 
Soient $A_D\in gl(n^2,\mathbb{C})$ diagonale et $A_N\in gl(n^2,\mathbb{C})$ nilpotente telles que

\begin{equation}(2i\pi\langle m,\omega\rangle-\mathcal{A}_{L,L'})=A_D-A_N
\end{equation}

Alors $A_N$ est l'opérateur

$$A_N:B\mapsto (\tilde{A} P^{\mathcal{L}}_L)_{\mathcal{N}}
B-B (P^{\mathcal{L}}_{L'}
\tilde{A})_{\mathcal{N}}$$ 

\noindent De plus,

\begin{equation}(2i\pi\langle m,\omega\rangle-\mathcal{A}_{L,L'})^{-1}=A_D^{-1}(I+A_NA_D^{-1}+\dots + (A_NA_D^{-1})^{n^2-1})
\end{equation}

\noindent Nous allons estimer $(2i\pi\langle m,\omega\rangle-\mathcal{A}_{L,L'})^{-1}$, pour $m\in \mathbb{Z}^d$ si $L=\bar{L}'$ et 
$m\in \frac{1}{2}\mathbb{Z}^d$ si $L\neq \bar{L}'$.
Chaque coefficient de $A_D^{-1}(A_NA_D^{-1})^{j-1}$ est de la forme $\frac{p}{q}$ avec $\mid p\mid \leq \mid \mid A_N\mid \mid ^{j-1}$ et $q=\beta_1\dots \beta_j$ où 
les $\beta_i$ sont des valeurs propres de $ 2i\pi\langle m,\omega\rangle-\mathcal{A}_{L,L'}$. Or

\begin{equation}\sigma(\mathcal{A}_{L,L'})=\{\alpha-\alpha'\ \mid \ \alpha\in \sigma(\tilde{A}_{\mid L}),\alpha'\in \sigma(\tilde{A}_{\mid L'}) \}
\end{equation}

\noindent 
et de plus, pour tous $\alpha\in \sigma(\tilde{A}_{\mid L}),\alpha'\in \sigma(\tilde{A}_{\mid L'})$,

\begin{equation}\mid \alpha-\alpha'- 2i\pi\langle m,\omega\rangle\mid \geq \frac{\kappa'}{\mid m\mid ^\tau}
\end{equation}

pour tout $m\in \mathbb{Z}^d$ si $L=\bar{L}'$ et tout $m\in \frac{1}{2}\mathbb{Z}^d$ si $L\neq \bar{L}'$
Ainsi,

\begin{equation}\begin{split} \mid \mid (2i\pi\langle m,\omega\rangle-\mathcal{A}_{L,L'})^{-1}\mid \mid
&\leq c_n (1+\mid \mid      \tilde{A} _{\mathcal{N}}   \mid \mid \ (\mid \mid  P^{\mathcal{L}}_L\mid \mid \ 
+\mid \mid P^{\mathcal{L}}_{L'}  \mid \mid ))^{n^2-1}\left(\frac{\mid m\mid^\tau}{\kappa'}\right)^{n^2-1}
\end{split}\end{equation}

\noindent où $c_n$ ne dépend que de $n$, 
et \eqref{invR} implique \eqref{estimationreG}.

\noindent L'estimation \eqref{estimationreG} et le lemme \ref{gevrey} impliquent alors que

\begin{equation}||P^{\mathcal{L}'}_L\tilde{X}
P^{\mathcal{L}'}_{L'}||_{r'}
\leq C'' 
\left( \frac{1+||\tilde{A}_{\mathcal{N}}||}{(r-r')\kappa'}\right)^{D\gamma}
||P^{\mathcal{L}'}_L
\tilde{F}
P^{\mathcal{L}'}_{L'}||_r
\end{equation}

\noindent où $C'',D$ ne dépendent que de $n,d,\tau$. Ainsi,





\begin{equation}\begin{split}
||\tilde{X}||_{r'}&\leq \sum_{L,L'}||P^{\mathcal{L}'}_L
\tilde{X}
P^{\mathcal{L}'}_{L'}||_r
\leq C''
\left(\frac{1+||\tilde{A}_{\mathcal{N}}||}{(r-r')\kappa'}\right)^{D\gamma}
\sum_{L,L'}||P^{\mathcal{L}'}_L
\tilde{F} P^{\mathcal{L}'}_{L'}||_r\\
\end{split}\end{equation}

\noindent d'où

\begin{equation}||\tilde{X}||_{r'}\leq C_3
\left(\frac{1+||\tilde{A}||}{(r-r')\kappa'}\right)^{D'\gamma}
||\tilde{F}||_r\end{equation}

\noindent où $D',C_3$ ne dépendent que de $n,d,\tau$. $\Box$

\section{Lemme inductif}

\subsection{Lemmes auxiliaires}
%


%



\noindent Rappelons un lemme démontré dans \cite{C2} qui sera utilisé pour itérer le lemme inductif sans faire intervenir une 
fonction de renormalisation à chaque étape, ce qui améliorera de beaucoup les estimations.

\begin{lem}\label{NRdurable} Soient 
\begin{itemize}
\item $\kappa'\in ]0,1[$, $C>0$,
\item $\tilde{F}\in \mathcal{G}$,  
\item $\tilde{\epsilon}=||\tilde{F}||$,
\item $\tilde{N}\in \mathbb{N}$,
\item $\tilde{A}\in \mathcal{G}$ avec un spectre $DC^{\tilde{N}}_\omega(\kappa',\tau)$.

\end{itemize}

\noindent Il existe une constante $c$ ne dépendant que de $n\tau $ 
telle que si $\tilde{\epsilon}$ est assez petit pour que

\begin{equation}\label{cond1}\tilde{\epsilon}\leq c\left(\frac{C^\tau \kappa'}{1+||\tilde{A}||}\right) ^{2n}
\end{equation} 

\noindent et

\begin{equation}\label{cond2} \tilde{N}
\leq  \frac{|\log \tilde{\epsilon}|^4}{C}
\end{equation}

\noindent alors 
$\tilde{A}+{\tilde{F}}$ a un spectre $DC^{\tilde{N}}_\omega(\frac{3\kappa'}{4},\tau)$.
\end{lem}

\rem Comme dans \cite{C2}, si $G$ est un groupe compact, alors le lemme \ref{NRdurable} est vrai en remplaçant \eqref{cond1} par une condition de petitesse qui ne dépend pas de $\tilde{A}$.

\bigskip
\noindent Le lemme qui suit sera utilisé pour garantir qu'un seul doublement de période est nécessaire.

\begin{lem}\label{bppG} Soit $A,A'\in gl(n,\mathbb{R})$ et $H:2\mathbb{T}^d\rightarrow gl(n,\mathbb{R})$. 
Supposons que $H$ a de bonnes propriétés de périodicité par rapport à une $A$-décomposition $\mathcal{L}$ et que 

\begin{equation}\label{parité}\forall L,L'\in\mathcal{L}, P_L^{\mathcal{L}}(A'-A)P^{\mathcal{L}}_{L'}\neq 0
\Rightarrow P_L^{\mathcal{L}}HP^{\mathcal{L}}_{L'}\in C^0(\mathbb{T}^d)
\end{equation}

\noindent Alors $H$ a de bonnes propriétés de périodicité par rapport à une $A'$-décomposition moins fine que $\mathcal{L}$.

\end{lem}

La démonstration se trouve dans \cite{C2} (lemme 4.2).

\begin{soulem}\label{fouGevrey}Soit $f\in C^{G,2}_r(2\mathbb{T}^d,gl(n,\mathbb{C}))$. Alors pour tout 

\begin{equation}N\geq  \frac{e^{2}d^{2}}{r^{4}}
\end{equation}

\noindent on a l'estimation

\begin{equation}\sum_{\mid m\mid >N}  \mid\mid \hat{f}(m)\mid  \mid
\leq C_d \mid \mid f\mid \mid_rN^De^{-\sqrt[4]{N}}
\end{equation}

\noindent où $C_d $ ne dépend que de $d$.

\end{soulem}

\dem Par définition de $\mid\mid f\mid\mid_r$, on a pour tout $\alpha\in \mathbb{N}^d$:

\begin{equation}\mid  \mid\widehat{ \partial^\alpha f} (m)\mid  \mid\leq \sup_\theta \mid \mid \partial^\alpha f(\theta)\mid \mid
\leq \frac{\alpha!^2}{r^{2\mid \alpha \mid }}\mid \mid f\mid \mid_r
\end{equation}

\noindent Or 

\begin{equation}\partial^\alpha f(\theta)\sim \sum_m \hat{f}(m) \partial^\alpha(e^{2i\pi \langle m,\theta\rangle})
\sim \sum_m \hat{f}(m) \prod_j( 2i\pi m_j )^{\alpha_j}.(e^{2i\pi \langle m,\theta\rangle})
\end{equation}

\noindent donc

\begin{equation} \widehat{ \partial^\alpha f} (m)
=\prod_j ( 2i\pi m_j )^{\alpha_j}\hat{f}(m)\end{equation}

\noindent et donc

\begin{equation}\label{alphagen}\mid \mid \hat{  f} (m)\mid \mid
\leq \frac{\alpha!^2}{\prod_j \mid 2i\pi m_j \mid^{\alpha_j} r^{2\mid \alpha \mid }}\mid \mid f\mid \mid_r
\end{equation}
%



%


%

\noindent Il existe forcément $1\leq j\leq d$ tel que $\mid m_j\mid \geq \frac{\mid m \mid}{d}$. 
Notons $M=\mid m\mid$ et supposons que $M\geq 1$. Réécrivons \eqref{alphagen} avec $\alpha=(\alpha_1,\dots, \alpha_d)$ où $\alpha_l=E(\sqrt[4]{M})$ si $l=j$ et $\alpha_l=0$ sinon. Alors

\begin{equation}\mid \mid \hat{  f} (m)\mid \mid
\leq \frac{(\sqrt[4]{M})!^{2}}{(\frac{2\pi r^2 M}{d})^{\sqrt[4]{M}}} \mid \mid f\mid \mid_r
\leq c \frac{(\sqrt[4]{M}+1)^{2(\sqrt[4]{M}+1)}e^{-2\sqrt[4]{M}}}{(\frac{2\pi r^2 M}{d})^{\sqrt[4]{M}}} \mid \mid f\mid \mid_r
\leq c \frac{(2\sqrt{M})^{\sqrt[4]{M}+1}e^{-2\sqrt[4]{M}}}{(\frac{2\pi r^2 M}{d})^{\sqrt[4]{M}}} \mid \mid f\mid \mid_r
\end{equation}

\noindent où $c$ est une constante numérique.
%



%

\noindent Ceci implique que

\begin{equation}\sum_{\mid m\mid >N} \mid \mid \hat{f}(m)\mid  \mid
\leq C_d\sum_{M>N}\frac{M^{d+\frac{1}{2}}e^{-2\sqrt[4]{M}}}{(\frac{\sqrt{M}\pi r^2 }{d})^{\sqrt[4]{M}}} \mid \mid f\mid \mid_r
\end{equation}

\noindent où $C_d$ ne dépend que de $d$, d'où

\begin{equation}\sum_{\mid m\mid >N} \mid \mid \hat{f}(m)\mid  \mid
\leq C_d \mid \mid f\mid \mid_r \int_N^{+\infty} \frac{t^{d+\frac{1}{2}}}{(\frac{\sqrt{t}}{d}r^2)^{\sqrt[4]{t}}}e^{-2\sqrt[4]{t}}dt
\end{equation}

\noindent et en effectuant le changement de variable $T=t^\frac{1}{4}$,

\begin{equation}\begin{split}\sum_{\mid m\mid >N} \mid \mid \hat{f}(m)\mid \mid 
&\leq C'_d \mid \mid f\mid \mid_r \int_{\sqrt[4]{N}}^{+\infty} \frac{T^{4d+2}}{(\frac{T^2}{d}r^2)^{T}}e^{-2T}T^3 dT\\
\end{split}\end{equation}

\noindent pour $C'_d$ ne dépendant que de $d$; finalement, pour certains $C''_d,D$ ne dépendant que de $d$,

\begin{equation}\begin{split}\sum_{\mid m\mid >N} \mid  \mid\hat{f}(m)\mid  \mid&
\leq C''_d \mid \mid f\mid \mid_r \frac{N^D}{\mid \log \left(\frac{d}{N^\frac{1}{2}r^2}\right)\mid ^D}\left(\frac{d}{N^\frac{1}{2}r^2}\right)^{\sqrt[4]{N}}\\
\end{split}\end{equation}

\noindent et comme, par hypothèse, $\frac{d}{N^\frac{1}{2}r^2}\leq e^{-1}$, alors

\begin{equation}\begin{split}\sum_{\mid m\mid >N} \mid \mid \hat{f}(m)\mid  \mid
&\leq C''_d \mid \mid f\mid \mid_r N^De^{-\sqrt[4]{N}} \ \Box
\end{split}\end{equation}

\begin{lem}\label{troncG}Soient $0<r\leq 1$, $f\in C^{G,2}_{r}(2\mathbb{T}^d,gl(n,\mathbb{C}))$, $N\in\mathbb{N}$ et $f^N$ la troncation de $f$ à l'ordre 
$N$. Supposons que 

\begin{equation}N\geq \frac{16e^{2}d^{2}}{r^{4}}
\end{equation}

\noindent 
Alors 

\begin{equation}||f-f^N||_{\frac{r}{2}}\leq C_d \mid \mid f\mid\mid_{r} N^De^{-\sqrt[4]{N}}
\end{equation}

\noindent où $C_d$ et $D$ ne dépendent que de $d$.
\end{lem}

\dem Par définition,

\begin{equation} ||f-f^N||_{\frac{r}{2}}= \sup_{\alpha,\theta} \frac{(\frac{r}{2})^{2\mid \alpha\mid}} {\alpha!^2} \mid \partial^\alpha(f-f^N)(\theta)\mid
\end{equation}

\noindent donc

\begin{equation} \begin{split} ||f-f^N||_{\frac{r}{2}} &= \sup_{\alpha,\theta} \frac{(\frac{r}{2})^{2\mid \alpha\mid}} {\alpha!^2} 
\mid \sum_{m>N} \hat{f}(m) \partial^\alpha(e^{2i\pi\langle m,\theta\rangle})\mid\\
&= \sup_{\alpha,\theta} \frac{(\frac{r}{2})^{2\mid \alpha\mid}} {\alpha!^2} 
\mid \sum_{m>N} \hat{f}(m)m^\alpha e^{2i\pi\langle m,\theta\rangle}\mid\\
&= \sup_{\alpha,\theta} \frac{(\frac{r}{2})^{2\mid \alpha\mid}} {\alpha!^2} 
\mid \sum_{m>N} \widehat{\partial^\alpha f}(m) e^{2i\pi\langle m,\theta\rangle}\mid
\end{split}\end{equation}
%


\noindent et par le sous-lemme \ref{fouGevrey}, il existe $C_d,D$ ne dépendant que de $d$ tels que

\begin{equation} \begin{split} ||f-f^N||_{\frac{r}{2}}\leq C_d \sup_{\alpha} \frac{(\frac{r}{2})^{2\mid \alpha\mid}} {\alpha!^2}
\mid \mid\partial^\alpha f\mid\mid_{\frac{r}{2}} N^De^{-\sqrt[4]{N}}
\end{split}\end{equation}

\noindent Or, d'après \cite{MS}, lemme A.2, pour tout $\alpha$,

\begin{equation}\frac{(\frac{r}{2})^{2\mid \alpha\mid}}{\alpha!^2}\mid \mid\partial^\alpha f\mid\mid_{\frac{r}{2}} \leq \mid \mid f\mid\mid_{r}
\end{equation}

\noindent d'où

\begin{equation} \begin{split} ||f-f^N||_{\frac{r}{2}}\leq C_d \mid \mid f\mid\mid_{r}N^D e^{-\sqrt[4]{N}} \ \Box
\end{split}\end{equation}

\subsection{Premier lemme inductif}


\begin{prop}\label{iterG} 
Soient
\begin{itemize}
\item $\tilde{\epsilon}>0,\tilde{r}\leq 1$, $ \kappa'>0,\tilde{N}\in\mathbb{N},\gamma\geq n(n+1)$; 

\item ${\tilde{F}}\in C^\omega_{\tilde{r}}(2\mathbb{T}^d,\mathcal{G}), \tilde{A}\in \mathcal{G}$,
\item $\mathcal{L}$ une $(\tilde{A},\kappa',\gamma)$-décomposition. 
\end{itemize}
Il existe une constante $C''>0$ ne dépendant que de $\tau,n$ telle que si

\begin{enumerate}

\item $\tilde{A}$ a un spectre $DC^{\tilde{N}}_\omega(\kappa',\tau)$;

\item 

\begin{equation}\label{moyptG}||\hat{\tilde{F}}(0)||\leq  \tilde{\epsilon}\leq C''\left(\frac{\kappa'}{1+||\tilde{A}||}\right)^{2n}\end{equation}

\noindent et

\begin{equation}\label{cond2'G} (\frac{4ed}{\tilde{r}^2})^2\leq \tilde{N}
\leq  (4ed)^2|\log \tilde{\epsilon}|^4
\end{equation} 

\item $\tilde{F}$ a de bonnes propriétés de périodicité par rapport à $\mathcal{L}$

\end{enumerate}
 
\noindent alors il existe 
\begin{itemize}
\item $C'\in\mathbb{R}$ ne dépendant que de $n,d,\kappa,\tau$, 
\item $D\in\mathbb{N}$ ne dépendant que de $n,d,\tau$,
\item ${X}\in C^\omega_{\frac{\tilde{r}}{2}}(2\mathbb{T}^d, \mathcal{G})$, 
\item ${A}'\in \mathcal{G}$ 
\item une $(A',\frac{3\kappa'}{4 },\gamma) $-décomposition $\mathcal{L}'$
\end{itemize}
vérifiant les propriétés suivantes: 

\begin{enumerate}

\item \label{nrG} $A'$ a un spectre $DC^{\tilde{N}}_\omega(\frac{3\kappa'}{4},\tau)$,

\item \label{0-G} $||{A}'-\tilde{A}||\leq  \tilde{\epsilon}$;

\item \label{bpper-G} la fonction 
$F'\in C^\omega_{\frac{\tilde{r}}{2}}(2\mathbb{T}^d,\mathcal{G})$ définie par

\begin{equation}\label{4-G}
\forall \theta\in 2\mathbb{T}^d,\ \partial_\omega e^{{X}(\theta)}=(\tilde{A}+\tilde{F}(\theta))e^{{X}(\theta)}
-e^{{X}(\theta)}({A}'+{F}'(\theta))\end{equation}

\noindent a de bonnes propriétés de périodicité par rapport à $\mathcal{L}'$

\item \label{9-G} pour tout $s\leq \tilde{r}$,

\begin{equation}\label{pxpG}||X||_s
\leq C'\left(\frac{1+||\tilde{A}_{\mathcal{N}}||}{\kappa'(\tilde{r}-s)}\right)^{D\gamma}
||\tilde{F}||_{\tilde{r}}\end{equation}

\item \label{10-G} et pour tout $s\leq \frac{\tilde{r}}{2}$,


\begin{equation}\begin{split} ||F'||_s&\leq C'
e^{||X||_s}||\tilde{F}||_{\tilde{r}}(\tilde{N}^de^{- \sqrt[4]{\tilde{N}}}+\frac{( 1+||\tilde{A}_{\mathcal{N}}||)^{D\gamma}}{(\kappa'(\tilde{r}-s))^{D\gamma}}||\tilde{F}||_{\tilde{r}}
(1+e^{||X||_s}))\end{split}\end{equation}

\end{enumerate}

\noindent De plus, si $\tilde{F}$ est continue sur $\mathbb{T}^d$, alors $X$ et $F'$ le sont également.

\end{prop}

\dem $\bullet$ Par hypothèse,  
$\tilde{F}$ a de bonnes propriétés de périodicité par rapport à $\mathcal{L}$ et une certaine famille $(m_L)$ et $\tilde{A}$ a un spectre $DC^{\tilde{N}}_\omega(\kappa',\tau)$, donc on peut appliquer 
la proposition \ref{homolG}. Soit ${X}\in C^\omega_{r'}(2\mathbb{T}^d,\mathcal{G})$ une solution 
de

\begin{equation}\forall \theta\in 2\mathbb{T}^d,\ 
\partial_\omega {X}(\theta)=[\tilde{A},{X}(\theta)]+\tilde{F}^{\tilde{N}}(\theta)
-\hat{\tilde{F}}(0)
\end{equation}

\noindent vérifiant la conclusion de la proposition \ref{homolG}.


\bigskip
\noindent Posons ${A}':=\tilde{A}+\hat{\tilde{F}}(0)$. 

$\bullet$ On a bien $||\tilde{A}-A'||=||\hat{\tilde{F}}(0)||$, d'où la propriété \ref{0-G}. 




\noindent 
$\bullet$ De plus, soit $c$ la constante donnée par le lemme \ref{NRdurable}, 
et supposons que $C''\leq c $. 
Les hypothèses \eqref{moyptG} et \eqref{cond2'G} permettent d'appliquer 
le lemme \ref{NRdurable} avec $C=\frac{1}{(4ed)^2}$ pour déduire que $A'$ a un spectre $DC^{\tilde{N}}_\omega
(\frac{3\kappa'}{4},\tau)$ d'où la propriété \ref{nrG}.

\noindent $\bullet$ Soit $F'\in C^\omega_{r'}(2\mathbb{T}^d,\mathcal{G})$ la fonction définie par \eqref{4-G}. Alors

\begin{equation}\label{F'bar-}\begin{split}&{F}'=e^{-{X}}
(\tilde{F} -\tilde{F}^{\tilde{N}} )
+e^{-{X} }
\tilde{F} 
(e^{{X} }-Id)+(e^{-{X} }-Id)\hat{\tilde{F}}(0)
-e^{-{X} }\sum_{k\geq 2}\frac{1}{k!}\sum_{l=0}^{k-1}{X} ^l
(\tilde{F}^{\tilde{N}} -\hat{\tilde{F}}(0))
{X} ^{k-1-l}\end{split}\end{equation}

\noindent Nous allons appliquer le lemme \ref{bppG} avec $A=\tilde{A}$ et $G=F'$, pour obtenir la propriété \ref{bpper-G}. 
La fonction ${F}'$ 
a de bonnes propriétés de périodicité par rapport à $\mathcal{L}$ 
et une certaine famille $(m_L)$ 
car ${X}$ et 
$\tilde{F}$ les ont. De plus, comme ${\tilde{F}}$ a de bonnes propriétés de périodicité par rapport à 
$\mathcal{L}$, 

\begin{equation}P^{\mathcal{L}}_L\hat{\tilde{F}}(0)P^{\mathcal{L}}_{L'}\neq 0\Rightarrow 
P^{\mathcal{L}}_L\tilde{F}P^{\mathcal{L}}_{L'}\in C^0(\mathbb{T}^d)
\end{equation}

\noindent or 

\begin{equation}
P^{\mathcal{L}}_L\tilde{F}P^{\mathcal{L}}_{L'}\in C^0(\mathbb{T}^d)
\Rightarrow m_L-m_{L'} \in \mathbb{Z}^d
\Rightarrow P^{\mathcal{L}}_L{F}'P^{\mathcal{L}}_{L'}\in C^0(\mathbb{T}^d)
\end{equation}

%


\noindent donc l'hypothèse \eqref{parité} du lemme \ref{bppG} est vérifiée. D'après le lemme \ref{bppG}, ${F}'$ a donc de bonnes 
propriétés de périodicité par rapport à une $A'$-décomposition $\mathcal{L}'$ qui est moins fine que $\mathcal{L}$, donc $\mathcal{L}'$ est une $(\tilde{A},
\kappa',\gamma)$-décomposition. Comme c'est une $(\tilde{A},
\kappa',\gamma)$-décomposition, chaque sous-espace $L\in\mathcal{L}'$ vérifie 

\begin{equation} \mid \mid P^{\mathcal{L}'}_L\mid \mid 
\leq C_0 \left(\frac{1+\mid \mid \tilde{A}_{\mathcal{N}} \mid \mid }{\kappa'}\right)^{\gamma} 
\end{equation}

\noindent donc

\begin{equation} \mid \mid P^{\mathcal{L}'}_L\mid \mid 
\leq C_0 \left(\frac{1+\mid \mid A' _{\mathcal{N}}\mid \mid +2\tilde{\epsilon}}{\kappa'}\right)^{\gamma} 
\leq C_0 \left(\frac{1+\mid \mid {A}' _{\mathcal{N}}\mid \mid }{\frac{3\kappa'}{4}}\right)^{\gamma} 
\end{equation}

\noindent donc $\mathcal{L}'$ est une $(A',\frac{3\kappa'}{4},\gamma)$-décomposition, 
d'où la propriété \ref{bpper-G}.

\bigskip
\noindent $\bullet$ La propriété \ref{9-G} est donnée par la proposition \ref{homolG}.

\bigskip
\noindent $\bullet$
L'estimation de $||F'||_s$ vient de l'expression \eqref{F'bar-} et du lemme \ref{troncG}, que l'on peut appliquer avec $r=\tilde{r}, f=\tilde{F}, N=\tilde{N}$, grâce à l'hypothèse 
\eqref{cond2'G} et la propriété \ref{9-G}. $\Box$

\subsection{Deuxième lemme inductif}

\bigskip
\noindent 
On peut obtenir directement l'étape inductive en itérant deux fois le lemme inductif sans renormalisation. 
Posons

\begin{equation}\label{N'G}N(\epsilon)=(4ed)^2|\log \epsilon|^4\end{equation}

\noindent et 

\begin{equation}\kappa''(\epsilon)=\frac{\kappa}{(9n{N(\epsilon)})^{\tau}}\end{equation}

\noindent Faisons une hypothèse supplémentaire sur $\gamma$. Soit $\bar{\gamma}\geq n(n+1)$ ne dépendant que de $n$ tel que $C_0^{\frac{1}{2\bar{\gamma}}}\leq 2$.

\begin{prop}\label{iter2G} 
Soient 
\begin{itemize}
\item $A\in \mathcal{G}$,

\item $r\leq 1,\gamma\geq \bar{\gamma}$,

\item $\bar{A},\bar{F}\in C^\omega_{r}(2\mathbb{T}^d,\mathcal{G})$, $\Psi\in C^\omega_{r}(2\mathbb{T}^d,G)$, 
\item $|\bar{F}|_{r}= \epsilon$, 


\end{itemize}

\bigskip
\noindent Il existe $\tilde{C}'>0$ ne dépendant que de $n,d,\kappa,\tau$ et $D_1\in\mathbb{N}$ ne dépendant que de $n,d,\tau$ 
tels que si en notant $r'=\frac{1}{\mid \log \epsilon\mid ^2}$,

\begin{enumerate}

\item

\begin{equation}\label{epsAG}\epsilon 
\leq \tilde{C}'\left(\frac{\kappa''(\epsilon)r'}{2(||A||+1)}\right)^{D_1\gamma}
\end{equation}

\noindent et 

\begin{equation}\label{petitesse2}\frac{4}{\mid \log \epsilon \mid ^2}\leq r
\end{equation}

\item $\bar{A}$ est réductible à $A$ 
par $\Psi$, 

\item \label{bp-G}$\Psi^{-1}\bar{F}\Psi$ a de bonnes propriétés de périodicité par rapport à une $(A,\kappa''(\epsilon),\gamma)$-décomposition $\mathcal{L}$,

\item pour tout $s\leq r$, $|\Psi|_{s}\leq (\frac{1}{\epsilon})^{\frac{1}{96}}$ et 
$|\Psi^{-1}|_{s}\leq (\frac{1}{\epsilon})
^{\frac{1}{96}}$,

\end{enumerate}
 
\noindent alors il existe 
\begin{itemize}
\item $Z'\in C^\omega_{r'}(2\mathbb{T}^d, G)$, 
\item $\bar{A}',\bar{F}'\in C^\omega_{r'}(2\mathbb{T}^d,\mathcal{G})$, 
\item $\Psi'\in C^\omega_{r'}(2\mathbb{T}^d, G)$, 
\item ${A}'\in \mathcal{G}$ 
\end{itemize}
\noindent vérifiant les propriétés suivantes: 

\begin{enumerate}

\item \label{55G} $\bar{A}'$ est 
réductible par $\Psi'$ à ${A}'$, 

\item \label{bpperG} la fonction $(\Psi')^{-1}\bar{F}'\Psi'$ a de bonnes propriétés de périodicité par rapport à une $(A',\kappa''(\epsilon^\frac{5}{2}),2\gamma)$-décomposition $\mathcal{L}''$

\item \label{renG} pour tout $s\leq r'$, $|\Psi'|_{s}\leq (\frac{1}{\epsilon^\frac{5}{2}})^{\frac{1}{96}}$ et 
$|(\Psi')^{-1}|_{s}\leq 
(\frac{1}{\epsilon^\frac{5}{2}})^{\frac{1}{96}}$,


\item \label{4G}

\begin{equation}\partial_\omega Z' =(\bar{A} +\bar{F} )Z' 
-Z' (\bar{A}' +\bar{F}')\end{equation}

\item \label{00G} $||A'||\leq ||A||+\epsilon^{\frac{11}{12}}+8\pi {N}$;

\item \label{66G} 

\begin{equation}|Z'-Id|_{r'}\leq \frac{1}{\tilde{C}'}\left(\frac{2(1+||A||)|\log\epsilon|}
{r-r'}\right)^{D_1\gamma}\epsilon^{\frac{5}{6}}
\end{equation}

\noindent et

\begin{equation} 
|(Z')^{-1}-Id|_{r'}\leq \frac{1}{\tilde{C}'}\left(\frac{2(1+||A||)|\log\epsilon|}
{r-r'}\right)^{D_1\gamma}\epsilon^{\frac{5}{6}}
\end{equation}


\item \label{8G} $|\bar{F}'|_{r'}\leq \epsilon^\frac{5}{2}$,

\item \label{phitrivG} la fonction $\Psi'^{-1}\Psi$ est triviale par rapport à $\mathcal{L}_{A,\kappa''}$.
%



\end{enumerate}

\noindent De plus, en dimension 2, si $\bar{A},\bar{F}$ sont continus sur $\mathbb{T}^d$, 
que l'hypothèse \ref{bp-G} est remplacée par 

\bigskip
\noindent \ref{bp-G}' pour toute fonction $H$ continue sur 
$\mathbb{T}^d$, $\Psi H \Psi^{-1}$ est continue sur $\mathbb{T}^d$

\bigskip
\noindent alors 
$Z',\bar{A}',\bar{F}'$ sont continus sur $\mathbb{T}^d$ et la propriété \ref{bpperG} est remplacée par 

\bigskip
\noindent \ref{bpperG}' pour toute fonction $H$ continue sur 
$\mathbb{T}^d$, $\Psi' H (\Psi')^{-1}$ est continue sur $\mathbb{T}^d$.

\bigskip
\noindent Enfin, si $\mathcal{G}=gl(n,\mathbb{C})$ ou $u(n)$ et si $\bar{A},\bar{F},\Psi$ sont continus sur $\mathbb{T}^d$, alors 
$Z',\bar{A}',\bar{F}',\Psi'$ sont continus sur $\mathbb{T}^d$.
 
\end{prop}

\dem $\bullet$
Vérifions d'abord que l'on peut appliquer la proposition  \ref{iterG} avec 

\begin{equation}
\tilde{\epsilon}=\epsilon;\ 
\tilde{r}=r;\ 
\kappa'=\kappa'';\ 
\tilde{N}=N;\
\end{equation}

Soit $\Phi$ une renormalisation de $A$ d'ordre $N$ et $\tilde{A}\in \mathcal{G}$ 
telle  que 

\begin{equation}\forall \theta\in 2\mathbb{T}^d,\ \partial_\omega \Phi(\theta)=A\Phi(\theta)-\Phi(\theta)\tilde{A}
\end{equation}

\noindent Posons $\Psi'=\Psi \Phi$ et soit

$$\tilde{F}:=(\Psi')^{-1}\bar{F}\Psi'$$
 
\noindent $\bullet$ La matrice $\tilde{A}$ a donc un spectre $DC_\omega^{N}(\kappa'',\tau)$ et la première hypothèse de \ref{iterG} est vérifiée. 

\noindent $\bullet$
Par hypothèse, $\Psi^{-1}\bar{F}\Psi$ a de bonnes propriétés de périodicité par rapport 
à une $(A,\kappa'',\gamma)$-décomposition $\mathcal{L}$ et une certaine famille $(m_L)$. 
Par ailleurs $\Phi$ est triviale par rapport à $\mathcal{L}_{A,\kappa''}$. Comme $\mathcal{L} $ et $\mathcal{L}_{A,\kappa''}$ sont des $A$-décompositions, on peut définir une $A$-décomposition $\bar{\mathcal{L}}$ de la manière suivante: 

\begin{equation}L\in \bar{\mathcal{L}} \Leftrightarrow \exists L_1\in \mathcal{L},L_2\in \mathcal{L}_{A,\kappa''}\ \mid \ 
L=L_1\cap L_2
\end{equation}

\noindent 
$\bar{\mathcal{L}}$ est une $(A,\frac{\kappa''}{C_0},2\gamma)$-décomposition car 
$\mathcal{L} $ et $\mathcal{L}_{A,\kappa''}$ sont des $(A,\kappa'',\gamma)$-décompositions et donc

\begin{equation}\mid \mid P^{\bar{\mathcal{L}}}_L\mid \mid 
=\mid \mid P^{\mathcal{L}}_{L_1} P^{\mathcal{L}_{A,\kappa''}}_{L_2}\mid \mid 
\leq C_0^2 \left( \frac{1+\mid \mid A_{\mathcal{N}}\mid \mid }{\kappa''}\right)^{2\gamma}
\end{equation}

\noindent 
De plus, 
$\tilde{F}$ a de bonnes propriétés de périodicité par rapport à $\bar{\mathcal{L}}$. Comme $\bar{\mathcal{L}}$ est une $(A,\frac{\kappa''}{C_0},2\gamma)$-décomposition, c'est aussi une 
$(\tilde{A},\frac{\kappa''}{C_0},2\gamma)$-décomposition (puisque $A$ et $\tilde{A}$ ont la même partie nilpotente) et donc la troisième hypothèse de \ref{iterG} est 
vérifiée.

\noindent De plus,

\begin{equation}||\hat{\tilde{F}}(0)||\leq |\tilde{F}|_{0}\leq |\Phi|_0|\Phi^{-1}|_0
|\Psi|_0|\Psi^{-1}|_0|\bar{F}|_0 \end{equation}

\noindent Or d'après \eqref{normephiG}, pour tout $s'\geq 0$,

\begin{equation}\mid\mid  \Phi \mid \mid_{s'} \leq nCC_0\left(\frac{1+||A_{\mathcal{N}}||}{\kappa''}\right)^{n(n+1)}e^{4\pi Ns'^2}
\end{equation}
%

%


\noindent où $C$ ne dépend que de $n,d$. L'hypothèse \eqref{epsAG}, avec  $\tilde{C}'$ assez petit et $D_1$ assez grand, et le fait que $\mid \mid A_{\mathcal{N}}\mid \mid \leq \mid \mid A\mid \mid$, 
impliquent alors que pour tout $s'\leq r$,

\begin{equation}\mid\mid\Phi\mid\mid_{s'}\leq \epsilon^{-\frac{1}{96}}e^{4\pi {N}s'^2}\end{equation}

\noindent donc pour $r'= \frac{1}{\mid \log \epsilon \mid^2}$,

\begin{equation}\mid\mid\Phi\mid\mid_{3r'}\leq \epsilon^{-\frac{1}{96}}e^{36\pi {N}r'^2}\end{equation}

\noindent et donc

\begin{equation}\mid\mid\Psi\Phi\mid\mid_{3r'}
\leq \mid\mid\Psi\mid\mid_{3r'} \mid\mid\Phi\mid\mid_{3r'}\leq \epsilon^{-\frac{1}{48}} e^{36\pi r'^2{N}}
\leq c\epsilon^{-\frac{1}{48}} \leq \epsilon^{-\frac{1}{40}}
\end{equation}

\noindent où $c$ ne dépend que de $d$, et de même pour $\Phi^{-1}$ (d'où la propriété \ref{renG}). Donc

\begin{equation}||\hat{\tilde{F}}(0)|| \leq \epsilon^{1-{2(r-r')}}C_0^2\left(\frac{1+||A_{\mathcal{N}}||}{\kappa''}\right)^{2n(n+1)}\end{equation}

\noindent d'où

\begin{equation}||\hat{\tilde{F}}(0)|| \leq \epsilon^{1-2(r-r')-\frac{1}{48}}\end{equation}



\noindent
Soit $C',D$ donnés par la proposition \ref{iterG} (ne dépendant donc que de $n,d$ et de $\tau$).
L'hypothèse \eqref{epsAG}, qui implique \eqref{moyptG} avec 

\begin{equation}\tilde{C}'\leq C'^4,\ 
D_1\gamma \geq 64n(n(n-1)+2)\tau\end{equation}

\noindent et l'expression \eqref{N'G} qui implique \eqref{cond2'G},  
%
permettent donc de satisfaire la deuxième hypothèse de la proposition \ref{iterG}.  

\bigskip
\noindent $\bullet$ On obtient ainsi des fonctions 
$X\in C^\omega_{r'}(2\mathbb{T}^d,\mathcal{G})$, 
${F}'\in C^\omega_{r'}(2\mathbb{T}^d,G)$, et une matrice ${A}'\in \mathcal{G}$ 
telles que 

\begin{itemize}

\item $A'$ a un spectre $DC^{N}_\omega(\frac{3}{4}\left(\frac{\kappa''}{C_0 }\right),\tau)$,

\item 

\begin{equation}||A'-\tilde{A}||\leq  \epsilon^{\frac{23}{24}}\end{equation}
%

%


\item $\partial_\omega e^X=(\tilde{A}+\tilde{F})e^X
-e^X({A}'+{F}')
$,

\item ${F}'$ a de bonnes propriétés de périodicité par 
rapport à une $(A',\frac{3\kappa''}{4C_0 },2\gamma)$-décomposition $\mathcal{L}'$

\item

\noindent pour tout $s\leq r$, 

\begin{equation}\label{pxp}
| X|_{s}\leq C'\left(\frac{C_0(1+||A_{\mathcal{N}}||)}{\kappa''(r-s)}\right)^{D\gamma}
|\tilde{F}|_r
\end{equation}

\noindent et pour tout $s\leq \frac{r}{2}$,

\begin{equation}\label{pfp}\begin{split}|F'|_{s}&
\leq C'e^{|X|_{s}}|\tilde{F}|_{r}
(N^de^{-\sqrt[4]{N}}
+\left(\frac{2C_0(1+||A_{\mathcal{N}}||)}{\kappa''s}\right)^{D\gamma}|\tilde{F}|_{r}(1+e^{| X |_{s}}))\end{split}\end{equation}

\end{itemize}


\bigskip
\noindent Par hypothèse, $2r'\leq \frac{r}{2}$ donc $X,F'\in C^\omega_{2r'}(2\mathbb{T}^d,\mathcal{G})$. 


%

%

%

%

%


\bigskip
\noindent $\bullet$ Montrons maintenant les estimations qui nous permettrons de réitérer la proposition \ref{iterG}. D'après \eqref{pfp}, où l'on a pris $s=2r'\leq \frac{r}{2}$,

\begin{equation}\begin{split}|F'|_{2r'}&
\leq C'e^{|X|_{2r'}}|\tilde{F}|_{3r'}
(N^de^{-\sqrt[4]{N}}
+\left(\frac{C_0(1+||A_{\mathcal{N}}||)}{\kappa''r'}\right)^{D\gamma}|\tilde{F}|_{3r'}(1+e^{| X |_{2r'}}))\end{split}\end{equation}

\noindent Or

\begin{equation}|\tilde{F}|_{3r'}\leq |\Psi'|_{3r'}|\Psi'^{-1}|_{3r'}|\bar{F}|_r
\leq \epsilon^{1-\frac{1}{20}}\end{equation}

\noindent et d'après \eqref{pxp} avec $s=2r'$,

\begin{equation}| X|_{2r'}
\leq C'\left(\frac{C_0(1+||{A}_{\mathcal{N}}||)}{\kappa''r'}\right)^{D\gamma}
| \tilde{F}|_{3r'}\tag{\ref {pxp}}\end{equation}

\noindent donc d'après \eqref{epsAG}, si $\tilde{C}'$ est assez petit et $D_1$ assez grand en fonction de $n,D, C'$,
alors

\begin{equation}| X|_{2r'}\leq \epsilon^{\frac{5}{6}}\end{equation}

\noindent donc

\begin{equation}e^{| X|_{2r'}}\leq 2\end{equation}

%

%

\noindent donc par l'hypothèse \eqref{epsAG},

\begin{equation}|{F}'|_{2r'}\leq 2C' \epsilon^{\frac{19}{20}}(N^d\epsilon
+3\epsilon^{\frac{19}{20}})
\end{equation}

\noindent Il existe une constante $c_d$ ne dépendant que de $d$ telle que si $\epsilon \leq c_d$, alors 

\begin{equation} \mid \log \epsilon \mid^d \leq \epsilon^{\frac{1}{96}}
\end{equation}

\noindent et dans ce cas,

\begin{equation}|{F}'|_{2r'}
\leq 
\epsilon^{\frac{9}{5}}
\end{equation}
%

%

%

\noindent $\bullet$ Nous pouvons donc appliquer à nouveau la proposition \ref{iterG} avec

\begin{equation}
\tilde{\epsilon}=\epsilon^{\frac{9}{5}};\ 
\tilde{r}=2r';\ 
\kappa'=\frac{3}{4}\kappa'';\ 
\tilde{N}=N;\
\tilde{A}=A';\ 
\tilde{F}=F';\ 
\mathcal{L}=\mathcal{L}'.
\end{equation}

\noindent Soient $A''\in\mathcal{G},F''\in C^\omega_{r'}(2\mathbb{T}^d,\mathcal{G}),\mathcal{L}'', X''\in C^\omega_{r'}(2\mathbb{T}^d,\mathcal{G})$ tels que

\begin{itemize}

\item $A''$ a un spectre $DC^N_\omega(\frac{9}{32}\kappa'',\tau$,

\item $\mid \mid A''-A'\mid \mid \leq \epsilon^{\frac{9}{5}}$, 

\item $\partial_\omega e^{X''}= (A'+F')e^{X''}-e^{X''}(A''+F'')$,

\item $F''$ a de bonnes propriétés de périodicité par rapport à une $(A'', \frac{9}{32}\kappa''(\epsilon), 2\gamma)$-décomposition $\mathcal{L}''$,

\item 

\begin{equation}| X''|_{r'}
\leq C'\left(\frac{3C_0(1+||{A}_{\mathcal{N}}||)}{\kappa'' r'}\right)^{D\gamma}
| F'|_{2r'}\leq \epsilon^{\frac{8}{5}}\end{equation}

\item et

\begin{equation}\label{f''}\begin{split}|F''|_{r'}&
\leq C'e^{|X''|_{r'}}|F'|_{2r'}
(N^de^{-\sqrt[4]{N}}
+\left(\frac{C_0(1+||A_{\mathcal{N}}||)}{\kappa''r'}\right)^{D\gamma}|F'|_{2r'}(1+e^{| X'' |_{r'}}))\\
\end{split}\end{equation}

\end{itemize}

\noindent Posons $Z'=\Psi'e^Xe^{X''}\Psi'^{-1},\bar{A}'=(\partial_\omega \Psi'+\Psi' A'')\Psi'^{-1}, \bar{F}'=\Psi' F''\Psi'^{-1}$. 

\bigskip
\noindent $\bullet$ \eqref{f''} implique que 

\begin{equation}
\begin{split}|F''|_{r'}
&\leq 2C' \epsilon^{\frac{9}{5}}(\epsilon^\frac{11}{12}+\epsilon^\frac{8}{5})\leq \epsilon^3
\end{split}\end{equation}

\noindent et finalement

\begin{equation}\begin{split}|\Psi'F''\Psi'^{-1}|_{r'}&
\leq \epsilon^{\frac{5}{2}}
\end{split}\end{equation}

\noindent d'où \ref{8G}.

\bigskip
\noindent $\bullet$ Vérifions que $\mathcal{L}'$ est une $(A'',\kappa''(\epsilon^\frac{5}{2}),2\gamma)$-décomposition. Il suffit de vérifier que 

\begin{equation}
\frac{9\kappa''(\epsilon)}{32}\geq \kappa''(\epsilon^\frac{5}{2})
\end{equation}

\noindent c'est-à-dire

\begin{equation}
\frac{9\kappa}{32(9n{(4ed)^2|\log \epsilon|^4})^{\tau}}\geq \kappa''(\epsilon^\frac{5}{2})=\frac{\kappa}{(9n{N(\epsilon^\frac{5}{2})})^{\tau}}
=\frac{\kappa}{(9n{(4ed)^2|\log \epsilon^\frac{5}{2}|^4})^{\tau}}
\end{equation}

\noindent ou encore

\begin{equation}
\frac{9}{32}\geq (\frac{2}{5})^{4\tau}
\end{equation}

\noindent ce qui est vrai puisque $\tau\geq 1$, d'où \ref{bpperG}.

\bigskip
\noindent $\bullet$ Nous allons estimer $|\Psi\Phi e^Xe^{X''}(\Psi\Phi) ^{-1}-Id|_{r'}$. On a 

\begin{equation}\mid e^X e^{X''}-Id\mid _{r'}=\mid e^X (e^{X''}-e^{-X})\mid _{r'}\leq 2(\mid e^X -Id\mid _{r'}+\mid  e^{X''}-Id\mid _{r'})
\end{equation}

\noindent et de plus

\begin{equation}\begin{split}| e^X -Id|_{r'}&
\leq C'\left(\frac{2C_0(1+||A_{\mathcal{N}}||)}{\kappa''(\epsilon)r'}\right)^{D\gamma}|\tilde{F}|_r\\
&\leq C''\left(\frac{2C_0(1+||A_{\mathcal{N}}||){N(\epsilon)}^{\tau }}{\kappa r'}\right)^{D\gamma}|\tilde{F}|_r
\end{split}\end{equation}

\noindent pour un certain $C''$ ne dépendant que de $n,d,\kappa,\tau$, et de même pour $| e^{X''} -Id|_{r'}$. Donc

\begin{equation}|\Psi\Phi e^Xe^{X''}(\Psi\Phi) ^{-1}-Id|_{r'}\leq 4C''\left(\frac{2C_0(1+||A_{\mathcal{N}}||){N(\epsilon)}^{\tau }}{\kappa r}
\right)^{D\gamma}
\epsilon^{\frac{9}{10}}
\end{equation}

\noindent c'est-à-dire

\begin{equation}\begin{split}|\Psi\Phi e^Xe^{X''}(\Psi\Phi) ^{-1}-Id|_{r'}&\leq C_3\left(\frac{2(1+||A_{\mathcal{N}}||)|\log\epsilon|}
{\kappa r}\right)^{D'_1\gamma}
\epsilon^{\frac{9}{10}}
\end{split}\end{equation}

\noindent pour un certain $C_3$ ne dépendant que de $n,d,\kappa,\tau$ et $D'_1$ ne dépendant que de $n,d,\tau$. 
On obtient la même estimation pour $|\Psi\Phi e^{-X''}e^{-X}(\Psi\Phi) ^{-1}-Id|_{r'}$, donc la propriété \ref{66G} est vérifiée avec $\tilde{C}'$ assez petit et $D_1$ assez grand en fonction de 
$n,d,\kappa,\tau$.


\bigskip
\noindent $\bullet$ De plus, 

\begin{equation}||A''||\leq ||A''-A'||+||A'-A||+||A||\leq ||A||+\epsilon^{\frac{11}{12}}+8\pi {N}
\end{equation}

\noindent d'où la propriété \ref{00G}.

\subsubsection{Cas de la dimension 2}

\noindent En dimension 2, comme par hypothèse $\Psi^{-1}\bar{F}\Psi$ est continue sur $\mathbb{T}^d$, et par construction de $\Phi$, 
alors $\tilde{F}$,${X}$ et ${F}'$ sont continues sur $\mathbb{T}^d$. 

\noindent Donc les fonctions $\Phi {F}'\Phi^{-1},\Phi\hat{\tilde{F}}(0)\Phi^{-1}$ et 
$\Phi {X}\Phi^{-1}$ sont continues 
sur $\mathbb{T}^d$, et par hypothèse sur 
$\Psi$, alors $\Psi\Phi {F}'(\Psi\Phi)^{-1},\Psi\Phi\hat{\tilde{F}}(0)(\Psi\Phi)^{-1}$ et $\Psi\Phi
 {X}(\Psi\Phi)^{-1}$ sont donc continues 
sur $\mathbb{T}^d$ et donc $\bar{A}'=\bar{A}+\Psi\Phi\hat{\tilde{F}}(0)(\Psi\Phi)^{-1}$ est continue sur $\mathbb{T}^d$.

\noindent Il ne reste qu'à vérifier que pour toute fonction $H$ continue sur $\mathbb{T}^d$, la fonction 
$(\Psi\Phi)^{-1}H\Psi\Phi$ est continue sur $\mathbb{T}^d$. Mais

$$(\Psi\Phi)^{-1}H\Psi\Phi=\Phi^{-1}\Psi^{-1}H\Psi\Phi$$ 

\noindent Par hypothèse, $\Psi^{-1}H\Psi$ est continue sur $\mathbb{T}^d$, et donc $\Phi^{-1}\Psi^{-1}H\Psi\Phi$ aussi. 
$\Box$

\section{Itération}

\begin{lem}\label{eps'1G}

\noindent Soit $C'\leq 1,b_0>0,\tau\geq 1$. Soit $D_2,\gamma_0\in\mathbb{N}$. Il existe $C$ ne dépendant que de $C',D_2,\gamma_0,\tau$
 tel que 
pour tout

\begin{equation}\label{lemnum}\epsilon\leq C\left(\frac{1}{b_0+1}\right)^{16\gamma_0D_2}
\end{equation}

\noindent et en posant pour tout $k$

\begin{equation}\left\{\begin{array}{c}
\epsilon_k=  \epsilon_{0}^{(\frac{5}{2})^k}\\
\gamma_k=2^k\gamma_0\\
r_k=\frac{1}{\mid \log \epsilon_{k-1}\mid^2}\\
b_k:=b_{k-1}+\frac{1}{C'}\mid \log \epsilon_{k-1}\mid^4\\
\kappa_k= \frac{C'}{ \mid \log \epsilon_k\mid ^{4\tau}}
\end{array}\right.\end{equation}

\noindent  
alors pour tout $k\in\mathbb{N}$,


%

\begin{equation}\label{epsk2}\frac{4}{r_k}\leq \mid \log \epsilon_k\mid^2
\end{equation}

\noindent et

\begin{equation}\label{epsk1}
a_k:=\left(\frac{b_k+1 }{\kappa_k r_{k+1}}\right)^{D_2\gamma_k}\epsilon_k\leq C'
\end{equation}


\end{lem}

\dem 
La propriété \eqref{epsk2} est évidente. 

%


\noindent \eqref{epsk1} est vraie si 

\begin{equation}\label{epsk1'}
\left(\frac{2(1+b_0)\mid \log \epsilon_0\mid^{10\tau} \sum_{j=0}^{k-1}(\frac{5}{2})^{j+2k} }{C'^2 }\right)^{2^kD_2\gamma_0}\epsilon_0^{(\frac{5}{2})^k}\leq C'
\end{equation}

\noindent et il suffit pour cela que

\begin{equation}\label{epsk1''}
\left(\frac{2(1+b_0)\mid \log \epsilon_0\mid^{10\tau} (\frac{5}{2})^{4k} }{C'^2 }\right)^{D_2\gamma_0}\epsilon_0^{(\frac{5}{4})^k}\leq C'
\end{equation}

\noindent ce qui est vrai pour tout $k$ si $\epsilon_0$ est assez petit pour que

\begin{equation}
(5(1+b_0))^{D_2\gamma_0}\mid \log \epsilon_0\mid ^{80\tau D_2\gamma_0}\epsilon_0\leq C'^{16D_2\gamma_0+1}\ \Box
\end{equation}

\begin{thm}\label{PRG} Soit $r\leq 1,A\in \mathcal{G}$ et $F\in C^{G,2}_r(2\mathbb{T}^d, \mathcal{G})$ avec de bonnes propriétés de périodicité par rapport à 
$\mathcal{L}_A$. 

\noindent 
Il existe $C$ ne dépendant que de $n,d,\tau,\kappa$ et $D_3$ ne dépendant que de $n,d,\tau,\kappa,A$
 tel que si 


\begin{equation}\label{petitesse3}|F|_r\leq\epsilon_0'= \left(\frac{C}{\mid \mid A\mid \mid +1}\right)^{D_3}
\end{equation}

\noindent alors 
pour tout $\epsilon\leq \epsilon_0'$, il existe 

\begin{itemize}
\item $r_\epsilon>0$,
\item $Z_\epsilon\in C^{G,2}_{r_\epsilon}(2\mathbb{T}^d,G)$,
\item $A_\epsilon\in \mathcal{G}$, 
\item $\bar{A}_\epsilon,\bar{F}_\epsilon\in C^{G,2}_{r_\epsilon}
(2\mathbb{T}^d,\mathcal{G})$, 
\end{itemize}
 
\noindent tels que 
\begin{enumerate}
\item \label{Ared} $\bar{A}_\epsilon$ est réductible à ${A}_\epsilon$,

\item \label{3'} $|\bar{F}_\epsilon|_{r_\epsilon}\leq \epsilon$

\item \label{2'} pour tout $\theta\in 2\mathbb{T}^d$,

$$\partial_\omega Z_\epsilon(\theta)=(A+F(\theta))Z_\epsilon(\theta)-Z_\epsilon(\theta)(\bar{A}_\epsilon(\theta)
+\bar{F}_\epsilon(\theta))
$$

\item \label{1'} 

$$|Z_\epsilon-Id|_{r_\epsilon}\leq 2^{D_3}\epsilon_0^{\frac{1}{4}-4r_\epsilon}$$ 

\noindent et 

$$|Z_\epsilon^{-1}-Id|_{r_\epsilon}\leq 2^{D_3}\epsilon_0^{\frac{1}{4}-4r_\epsilon}$$

\item \label{4''} $Z_\epsilon, \partial_\omega Z_\epsilon$ sont bornées dans $C^\omega_{r_\epsilon}(2\mathbb{T}^d,
gl(n,\mathbb{C}))$ indépendamment de $\epsilon$;

\end{enumerate}

\noindent De plus, en dimension 2 ou si $\mathcal{G}=gl(n,\mathbb{C})$ ou $u(n)$, si $F$ est continu sur $\mathbb{T}^d$, alors 
$\bar{A}_\epsilon, \bar{F}_\epsilon$ et $Z_\epsilon$ sont en faits continus sur $\mathbb{T}^d$. Si $\mathcal{G}=o(n)$ ou $u(n)$, alors $D_3$ est indépendant de $A$.

\end{thm}

\dem La preuve se fait par itération de la proposition \ref{iter2G}, grâce au lemme \ref{eps'1G}. Soit $\gamma\geq \bar{\gamma}$ ne dépendant que de $A, \kappa,n$ tel que $\mathcal{L}_A$ soit une 
$(A,\kappa,\gamma)$-décomposition (si $\mathcal{G}=o(n)$ ou $u(n)$, on peut prendre $\gamma=\bar{\gamma}$ qui est indépendant de $A$). Appliquons la proposition \ref{iter2G} avec $\Psi \equiv Id$. Soient $\tilde{C}', D_1$ donnés par la proposition \ref{iter2G}. Nous allons appliquer le lemme \ref{eps'1G} avec $C'=\tilde{C}', b_0=\mid \mid A\mid \mid, D_2=D_1,\gamma_0=\gamma$. Soit $C$ donné par le lemme \ref{eps'1G}, $D_3=16\gamma_0D_2$ et soit $\epsilon_0'$ 
satisfaisant la condition \eqref{petitesse3}, qui implique \eqref{lemnum}. 
Alors 
$\mathcal{L}_A$ est en particulier une $(A,\kappa''(\epsilon_0'),\gamma)$-décomposition et $\epsilon_0'$ satisfait \eqref{epsAG} et \eqref{petitesse2}. On peut donc itérer la proposition \ref{iter2G} un nombre $k_\epsilon$ de fois qui est tel que 

\begin{equation}\epsilon_0'^{(\frac{5}{2})^{k_\epsilon}}\leq \epsilon
\end{equation}

\noindent pour obtenir $r_\epsilon, Z_\epsilon,A_\epsilon, \bar{A}_\epsilon, \bar{F}_\epsilon$ satisfaisant les propriétés \ref{Ared}, \ref{3'} et \ref{2'}, et où $Z_\epsilon$ est un produit de transformations $Z_k$ telles que, pour tout $k\leq 
k_\epsilon$,

\begin{equation}|Z_k-Id|_{r_\epsilon}\leq \frac{1}{\tilde{C}'}\left(\frac{2(1+||A||)}{r}5^k|\log\epsilon_0'|\right)^{D_1\gamma}\epsilon_0'^{\frac{5}{6}(\frac{5}{2})^k}
\end{equation}

\noindent et de même pour $Z_k^{-1}$. Les estimations \ref{1'} s'obtiennent par une récurrence simple et la propriété \ref{4''} se démontre par  une estimation de Cauchy. $\Box$

\bigskip
\noindent Comme corollaire, on obtient la quasi-densité des cocycles réductibles au voisinage d'une constante, dans la topologie Gevrey.

\begin{cor}Soit $\mathcal{G}$ une algèbre de Lie parmi $gl(n,\mathbb{C}), gl(n,\mathbb{R}), sp(n,\mathbb{R}), sl(n,\mathbb{R}),o(n),u(n)$, soit 
$A\in \mathcal{G}$, $r\leq 1$, $F\in C^{G,2}_r(\mathbb{T}^d,\mathcal{G})$. Il existe $\epsilon_0$ ne dépendant que de $n,d,\kappa,\tau,A,r$ tel que si 
$\mid\mid F-A\mid\mid_r\leq \epsilon_0$, alors pour tout $\epsilon>0$ il existe $r_\epsilon>0$ et $H\in C^{G,2}_{r'}(2\mathbb{T}^d,\mathcal{G})$ réductible tel que 
$\mid\mid F-H\mid \mid_{r_\epsilon}\leq \epsilon$. Si $\mathcal{G}$ est complexe ou si $n=2$, alors on peut supposer $H$ défini sur $\mathbb{T}^d$. Si $\mathcal{G}=o(n)$ ou $u(n)$, alors $\epsilon_0$ ne dépend pas de $A$.

\end{cor}

Ceci prouve les théorèmes \ref{th1} et \ref{th3}.

\end{document}